\title{On the homology of almost Calabi-Yau algebras \\ associated to $SU(3)$ modular invariants}
\author{
        David E. Evans and Mathew Pugh \\ \\
        School of Mathematics, \\
        Cardiff University, \\
        Senghennydd Road, \\
        Cardiff, CF24 4AG, \\
        Wales, U.K.
}
\date{\today}

\documentclass[12pt]{article}
\usepackage{amssymb}
\usepackage{graphicx}
\usepackage[usenames]{color}

\newtheorem{Def}{Definition}[section]
\newtheorem{Prop}[Def]{Proposition}

\newtheorem{Thm}[Def]{Theorem}

\begin{document}
\maketitle

\begin{abstract}
We compute the Hochschild homology and cohomology, and cyclic homology, of almost Calabi-Yau algebras for $SU(3)$ $\mathcal{ADE}$ graphs. These almost Calabi-Yau algebras are a higher rank analogue of the pre-projective algebras for Dynkin diagrams, which are $SU(2)$-related constructions.
The Hochschild (co)homology and cyclic homology of $A$ can be regarded as invariants for the braided subfactors associated to the $SU(3)$ modular invariants.
\end{abstract}

\section{Introduction}

The classical McKay correspondence appears in various contexts. It relates finite subgroups $\Gamma$ of $SU(2)$ with the algebraic geometry of the quotient Kleinian singularities $\mathbb{C}^2/\Gamma$ \cite{reid:2002} but also with the classification of $SU(2)$ modular invariants
and quantum subgroups of $SU(2)$ \cite{cappelli/itzykson/zuber:1987ii, zuber:2002, ocneanu:2000ii, ocneanu:2002, xu:1998, bockenhauer/evans:1999i, bockenhauer/evans:1999ii, bockenhauer/evans/kawahigashi:1999, bockenhauer/evans/kawahigashi:2000}.

Minimal resolutions of Kleinian singularities can be described via the moduli space of representations of the preprojective algebra associated to the action of $\Gamma$ \cite{crawley-boevey/holland:1998}.
Preprojective algebras associated to graphs were introduced in \cite{gelfand/ponomarev:1979}, and it was shown that they are finite dimensional if and only if the graphs are of $ADE$ type, that is, one of the simply laced Dynkin diagrams.
The preprojective algebra $A$ for an $ADE$ Dynkin diagram is a Frobenius algebra, that is, there is a linear function $f:A \rightarrow \mathbb{C}$ such that $(x,y):=f(xy)$ is a non-degenerate bilinear form (this is equivalent to the statement that $A$ is isomorphic to its dual $\widehat{A} = \mathrm{Hom}(A,\mathbb{C})$ as left (or right) $A$-modules). There is an automorphism $\beta$ of $A$, called the Nakayama automorphism of $A$ (associated to $f$), such that $(x,y) = (y,\beta(x))$, which yields an $A$-$A$ bimodule isomorphism $\widehat{A} \rightarrow {}_1 A_{\beta}$ \cite{yamagata:1996}.
The Nakayama automorphism for each $ADE$ graph was determined in \cite{erdmann/snashall:1998i, erdmann/snashall:1998ii} (see also \cite{brenner/butler/king:2002}).
The preprojective algebra $A$ has a finite resolution as an $A$-$A$ bimodule, which was used by Erdmann and Snashall to determine the Hochschild cohomology $HH^{\bullet}(A)$ of $A$ for the $ADE$ graphs $A_n$, along with its ring structure \cite{erdmann/snashall:1998i}, and $HH^2(A)$ for the graphs $D_n$ \cite{erdmann/snashall:1998ii}
This finite resolution yields a projective resolution of $A$ as an $A$-$A$ bimodule, which was used by Etingof and Eu to determine the Hochschild homology and cohomology, and cyclic homology, of $A$ for all $ADE$ Dynkin diagrams \cite{etingof/eu:2007}, along with the ring structure of the Hochschild cohomology \cite{eu:2007i}. The Hochschild homology and cohomology, cyclic homology, and ring structure of the Hochschild cohomology, for the preprojective algebra for the tadpole graphs $T_n$ were obtained in \cite{eu:2007ii}.
The Hochschild homology and cohomology for the case of the affine Dynkin diagrams, which are the McKay graphs for the finite subgroups of $SU(2)$, were determined in \cite{crawley-boevey/etingof/ginzburg:2007}.

More generally, one tries to understand singularities via a noncommutative algebra $A$, often called a noncommutative resolution, whose centre corresponds to the coordinate ring of the singularity \cite{vandenBergh:2004}. The algebra should be finitely generated over its centre, and the desired favourable resolution is the moduli space of representations of $A$, whose category of finitely generated modules is derived equivalent to the category of coherent sheaves of the resolution.
In the case of a quotient singularity $\mathbb{C}^3/\Gamma$ for a finite subgroup $\Gamma$ of $SU(3)$, the corresponding noncommutative algebra $A$ is a Calabi-Yau algebra of dimension 3.

Calabi-Yau algebras arise naturally in the study of Calabi-Yau manifolds, providing a noncommutative version of conventional Calabi-Yau geometry.
An algebra $A$ is Calabi-Yau of dimension $n$ if the bounded derived category of the abelian category of finite dimensional $A$-modules is a Calabi-Yau category of dimension $n$. In this case the global dimension of $A$ is $n$ \cite{bocklandt:2008}.
The derived category of coherent sheaves over an $n$-dimensional Calabi-Yau manifold is a Calabi-Yau category of dimension $n$ and they appear naturally in the study of boundary conditions of the $B$-model in superstring theory over the manifold. For more on Calabi-Yau algebras, see e.g. \cite{bocklandt:2008, ginzburg:2006}.

In \cite[Remark 4.5.7]{ginzburg:2006} Ginzburg introduced, in his terminology, $q$-deformed Calabi-Yau algebras. In the case where $q$ is not a root of unity, these algebras are Calabi-Yau algebras of dimension 3.
We study these algebras in the case where $q$ is a root of unity, which are the $SU(3)$ generalizations of preprojective algebras for the Coxeter-Dynkin diagrams $ADE$. We call these algebras \emph{almost Calabi-Yau algebras}. In a recent work \cite{evans/pugh:2010ii}, we determined the Nakayama automorphism for each $\mathcal{ADE}$ graph, and constructed a finite resolution of $A$ as an $A$-$A$ bimodule, see (\ref{exact_seq-almostCY}).

Our interest in these almost Calabi-Yau algebras came from subfactor theory, and in particular, braided subfactors of von Neumann algebras, which provide a framework for studying two dimensional conformal field theories and their modular invariant partition functions.
In the case of Wess-Zumino-Witten models associated to $SU(n)$ at level $k$, the Verlinde algebra is a non-degenerately braided system of endomorphisms ${}_N \mathcal{X}_N$, labelled by the positive energy representations of the loop group of $SU(n)_k$ on a type $\mathrm{III}_1$ factor $N$, with fusion rules $\lambda \mu = \bigoplus_{\nu} N_{\lambda \nu}^{\mu} \nu$ which exactly match those of the positive energy representations \cite{wassermann:1998}. The fusion matrices $N_{\lambda} = [N_{\rho \lambda}^{\sigma}]_{\rho,\sigma}$ are a family of commuting normal matrices which give a representation themselves of the fusion rules of the positive energy representations of the loop group of $SU(n)_k$, $N_{\lambda} N_{\mu} = \sum_{\nu} N_{\lambda \nu}^{\mu} N_{\nu}$.
This family $\{ N_{\lambda} \}$ of fusion matrices can be simultaneously diagonalised:
$$N_{\lambda} = \sum_{\sigma} \frac{S_{\sigma, \lambda}}{S_{\sigma,0}} S_{\sigma} S_{\sigma}^{\ast},$$
where $0$ is the trivial representation, and the eigenvalues $S_{\sigma, \lambda}/S_{\sigma,0}$ and eigenvectors $S_{\sigma} = [S_{\sigma, \mu}]_{\mu}$ are described by the statistics $S$ matrix.
The key structure in the conformal field theory is the modular invariant partition function $Z$. In the subfactor setting this is realised by
a braided subfactor $N \subset M$ where trivial (or permutation) invariants in the ambient factor $M$ when restricted to $N$ yield $Z$. This would mean that the dual canonical endomorphism is in $\Sigma({}_N \mathcal{X}_N)$, i.e. decomposes as a finite linear combination of endomorphisms in ${}_N \mathcal{X}_N$.
Indeed if this is the case for the inclusion $N \subset M$, then the process of $\alpha$-induction allows us to analyse the modular invariant,
providing two extensions of $\lambda$ on $N$ to endomorphisms $\alpha^{\pm}_{\lambda}$ of $M$, such that the matrix $Z_{\lambda,\mu} = \langle \alpha_{\lambda}^+, \alpha_{\mu}^- \rangle$ is a modular invariant \cite{bockenhauer/evans/kawahigashi:1999, bockenhauer/evans:2000, evans:2003}.
The action of the system ${}_N \mathcal{X}_N$ on the $N$-$M$ sectors ${}_N \mathcal{X}_M$ produces a \emph{nimrep} (non-negative matrix integer representation of the fusion rules) $G_{\lambda} G_{\mu} = \sum_{\nu} N_{\lambda \nu}^{\mu} G_{\nu}$,
whose spectrum reproduces exactly the diagonal part of the modular invariant, i.e.
$$G_{\lambda} = \sum_{\sigma} \frac{S_{\sigma,\lambda}}{S_{\sigma,0}} \psi_{\sigma} \psi_{\sigma}^{\ast},$$
with the spectrum of $G_{\lambda}$ given by
$G_{\lambda} = \{ S_{\mu, \lambda}/S_{\mu,0} \textrm{ with multiplicity } Z_{\mu,\mu} \}$
\cite[Theorem 4.16]{bockenhauer/evans/kawahigashi:2000}.

The systems ${}_N \mathcal{X}_N$, ${}_N \mathcal{X}_M$, ${}_M \mathcal{X}_M$ are (the irreducible objects of) tensor categories of endomorphisms with the Hom-spaces as their morphisms. Thus ${}_N \mathcal{X}_N$ gives a braided modular tensor category, and ${}_N \mathcal{X}_M$ a module category.

In our work we have focused on braided subfactors associated to $SU(3)$ modular invariants, which are labeled by a family of graphs which we call the $SU(3)$ $\mathcal{ADE}$ graphs. The complete list of the $SU(3)$ $\mathcal{ADE}$ graphs are illustrated in \cite[Figures 5-9]{evans/pugh:2011}.
For positive integer $k < \infty$ we have a braided modular tensor category ${}_N \mathcal{X}_N = \{ \lambda_{(p,l)} | \; 0 \leq p,l,p+l \leq k \}$, a non-degenerately braided system of endomorphisms on a type $\mathrm{III}_1$ factor $N$, which is generated by $\rho = \lambda_{(1,0)}$ and its conjugate $\overline{\rho} = \lambda_{(0,1)}$, where the irreducible endomorphisms $\lambda_{(p,l)}$ satisfy the fusion rules of $SU(3)_k$:
\begin{equation} \label{eqn:fusion_ruleSU(3)}
\lambda_{(p,l)} \otimes \rho \cong \lambda_{(p,l-1)} \oplus \lambda_{(p-1,l+1)} \oplus \lambda_{(p+1,l)},
\quad
\lambda_{(p,l)} \otimes \overline{\rho} \cong \lambda_{(p-1,l)} \oplus \lambda_{(p+1,l-1)} \oplus \lambda_{(p,l+1)},
\end{equation}
where $\lambda_{(p',l')}$ is understood to be zero if $p'<0$, $l'<0$ or $p'+l' \geq k+1$.
Then a pair $(\mathcal{G},W)$, of a cell system $W$ (see Section \ref{sect:almostCYalg}) on an $SU(3)$ $\mathcal{ADE}$ graph $\mathcal{G}$ with Coxeter number $k+3$ yields a braided subfactor $N \subset M$ and a module category ${}_N \mathcal{X}_M$, where the associated modular invariant, labeled by $\mathcal{G}$, is at level $k$.
For such a braided subfactor, the almost Calabi-Yau algebra can be constructed via a monoidal functor $F$, which is essentially the module category ${}_N \mathcal{X}_M$, from the $A_2$-Temperley-Lieb category to the category $\mathrm{Fun}({}_N \mathcal{X}_M,{}_N \mathcal{X}_M)$ of additive functors from ${}_N \mathcal{X}_M$ to itself.

The $A_2$-Temperley-Lieb category constructed in \cite{evans/pugh:2010ii} used ideas from planar algebras, and in particular, the $A_2$-planar algebras of \cite{evans/pugh:2009iii} (see also an earlier construction of the $A_2$-Temperley-Lieb category in \cite{cooper:2007}, and of the Temperley-Lieb category in \cite{turaev:1994, yamagami:2003}).

\begin{figure}[bt]
\begin{center}
  \includegraphics[width=40mm]{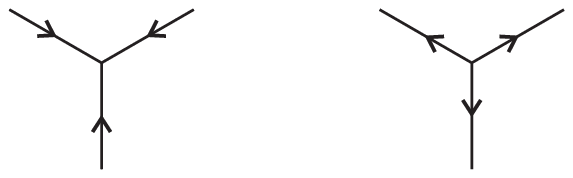}\\
 \caption{$A_2$ webs}\label{fig:A2-webs}
\end{center}
\end{figure}

For $m_i,n_i \geq 0$, an $A_2$-$(m_2,n_2),(m_1,n_1)$-tangle $T$ is a tangle on an rectangle with $m_2+n_2$, $m_1+n_1$ vertices along the top, bottom edges respectively, generated by $A_2$ webs (see Figure \ref{fig:A2-webs}) such that every free end of $T$ is attached to a vertex along the top or bottom of the rectangle in a way that respects the orientation of the strings, every vertex has a string attached to it, and the tangle contains no elliptic faces. We call a vertex a source vertex if the string attached to it has orientation away from the vertex. Similarly, a sink vertex will be a vertex where the string attached has orientation towards the vertex. Along the top, bottom edge the first $m_i$ vertices are source, sink vertices respectively, and the last $n_i$ are sink, source vertices respectively.
Let $V^{A_2}_{(m_2,n_2),(m_1,n_1)}$ be the quotient of the free vector space over $\mathbb{C}$ with basis the $A_2$-$(m_2,n_2),(m_1,n_1)$-tangles, by the Kuperberg ideal generated by the Kuperberg relations K1-K3 \cite{kuperberg:1996}.
Then at level $k$, the $A_2$-Temperley-Lieb category is defined to be a quotient of the category $A_2\textrm{-}TL = \mathrm{Mat}(C^{A_2})$ by the negligible morphisms, where $C^{A_2}$ is the tensor category whose objects are projections in $V^{A_2}_{(m,n),(m,n)}$ and whose morphisms are $\mathrm{Hom}(p_1,p_2) = p_2 V^{A_2}_{(m_2,n_2),(m_1,n_1)} p_1$, for projections $p_i \in V^{A_2}_{(m_i,n_i),(m_i,n_i)}$, $i=1,2$. We write $A_2\textrm{-}TL_{(m,n)} = V^{A_2}_{(m,n),(m,n)}$, and $\rho$, $\overline{\rho}$ for the identity projections in $A_2\textrm{-}TL_{(1,0)}$, $A_2\textrm{-}TL_{(0,1)}$ respectively consisting of a single string with orientation downwards, upwards respectively. Then the identity diagram in $A_2\textrm{-}TL_{(m,n)}$, given by $m+n$ vertical strings where the first $m$ strings have downwards orientation and the next $n$ have upwards orientation, is expressed as $\rho^m \overline{\rho}^n$. It is a linear combination of simple projections $f_{(i,j)}$ for $i,j \geq 0$, $0 \leq i+j < m+n$ such that $i-j \cong m-n \textrm{ mod } 3$, and a simple projection $f_{(m,n)}$, where $f_{(1,0)} = \rho$, $f_{(0,1)} = \overline{\rho}$ and $f_{(0,0)}$ is the empty diagram. The morphisms $\mathfrak{f}_{(p,l)} = \mathrm{id}_{f_{(p,l)}}$ are generalized Jones-Wenzl projections.
The $f_{(p,l)}$ satisfy the fusion rules for $SU(3)$ \cite{evans/pugh:2010ii}:
\begin{equation} \label{eqn:fusion_rule-f(k,l)}
f_{(p,l)} \otimes \rho \cong f_{(p,l-1)} \oplus f_{(p-1,l+1)} \oplus f_{(p+1,l)},
\quad
f_{(p,l)} \otimes \overline{\rho} \cong f_{(p-1,l)} \oplus f_{(p+1,l-1)} \oplus f_{(p,l+1)}.
\end{equation}

At level $k$, the negligible morphisms are the ideal $\langle \mathfrak{f}_{(p,l)} | p+l=k+1 \rangle$ generated by $\mathfrak{f}_{(p,l)}$ such that $p+l=k+1$.
The $A_2$-Temperley-Lieb category is the quotient $A_2\textrm{-}TL^{(k)} := A_2\textrm{-}TL/ \langle \mathfrak{f}_{(p,l)} | p+l=k+1 \rangle$, which is semisimple with simple objects $f_{(p,l)}$, $p,l \geq 0$ such that $p+l \leq k$ which satisfy the fusion rules (\ref{eqn:fusion_ruleSU(3)}) of $SU(3)_k$, that is we have (\ref{eqn:fusion_rule-f(k,l)}) where $f_{(p',l')}$ is understood to be zero if $p'<0$, $l'<0$ or $p'+l' \geq k+1$.
The $A_2$-Temperley-Lieb category $A_2\textrm{-}TL^{(k)}$ may be identified with the braided modular tensor category ${}_N \mathcal{X}_N$, where the object $f_{(p,l)} \in A_2\textrm{-}TL^{(k)}$ is identified with $\lambda_{(p,l)} \in {}_N \mathcal{X}_N$.

Then the monoidal functor $F$ is given on the simple objects $f_{(p,l)}$ of $A_2\textrm{-}TL^{(k)}$ by
\begin{equation} \label{eqn:functorF}
F(f_{(p,l)}) = \bigoplus_{i,j \in \mathcal{G}_0} G_{\lambda_{(p,l)}}(i,j) \, \mathbb{C}_{i,j},
\end{equation}
where
$\mathbb{C}_{i,j}$ are 1-dimensional $R$-$R$ bimodules, where $R = (\mathbb{C}\mathcal{G})_0$. The category of $R$-$R$ bimodules has a natural monoidal structure given by $\otimes_R$.
The functor $F$ is defined on the morphisms of $A_2\textrm{-}TL^{(k)}$ using the cell system $W$ and the Perron-Frobenius eigenvector of $\mathcal{G}$, see \cite[Section 2.9]{evans/pugh:2010ii}.

If $\mathcal{G}^{\mathrm{op}}$ denotes the opposite graph of $\mathcal{G}$ obtained by reversing the orientation of every edge of $\mathcal{G}$, we have that $F(\rho^m \overline{\rho}^n)$ is the $R$-$R$ bimodule with basis given by all paths of length $m+n$ on
$\mathcal{G}$, $\mathcal{G}^{\mathrm{op}}$, where the first $m$ edges are on $\mathcal{G}$ and the last $n$ edges are on $\mathcal{G}^{\mathrm{op}}$.
In particular $F(\rho^m) = (\mathbb{C}\mathcal{G})_m$, so that we have the graded algebra $\bigoplus_m F(\rho^m) = (\mathbb{C}\mathcal{G})$, the path algebra of $\mathcal{G}$.
The endomorphisms $\rho^m$ are not irreducible however, but decompose into direct sums of the generalized Jones-Wenzl projections $f_{(p,0)}$.
The natural algebra to consider is thus the graded algebra $\Sigma = \bigoplus_j F(f_{(j,0)})$, where the $p^{\mathrm{th}}$ graded part is $\Sigma_p = F(f_{(p,0)})$. The multiplication $\mu$ is defined by $\mu_{p,l} = F(\mathfrak{f}_{(p+l,0)}): \Sigma_p \otimes_R \Sigma_l \rightarrow \Sigma_{p+l}$, where $\mathfrak{f}_{(p,l)} = \mathrm{id}_{f_{(p,l)}}$.
The graded algebra $\Sigma$ is isomorphic to the almost Calabi-Yau algebra $A = A(\mathcal{G},W)$ \cite{cooper:2007, evans/pugh:2010ii}.

In Section \ref{sect:almostCYalg} we introduce the almost Calabi-Yau algebra $A = A(\mathcal{G},W)$ for a pair $(\mathcal{G},W)$ of a cell system $W$ on an $SU(3)$ $\mathcal{ADE}$ graph $\mathcal{G}$. Then in Section \ref{sect:resolution_almostCYalg} we determine a periodic projective resolution of $A$ as an $A$-$A$ bimodule, starting from the finite resolution of $A$ determined in \cite[Theorem 5.1]{evans/pugh:2010ii}, which will be used to determine the Hochschild (co)homology and cyclic homology of $A$ in Sections \ref{sect:Hoch_hom}-\ref{sect:Hoch_cohom}.
In Section \ref{sect:Hoch_hom-complex} we use the projective resolution determined in Section \ref{sect:resolution_almostCYalg} to construct a Hochschild homology complex for $A$, and introduce the cyclic homology of $A$ in Section \ref{sect:cyclic_hom}. We then determine the Hochschild and cyclic homology of $A$ in Sections \ref{sect:HH_0}-\ref{Sect:HH(A)-non-trivial_beta} for the graphs $\mathcal{A}^{(n)}$, $n = 4,5,6,7$, $\mathcal{D}^{(3k+3)}$, $k \geq 1$, $\mathcal{A}^{(n)\ast}$, $n \geq 5$, $\mathcal{D}^{(3k)\ast}$, $k \geq 2$, $\mathcal{E}^{(8)}$ and $\mathcal{E}^{(8)\ast}$.
Finally in Section \ref{sect:Hoch_cohom} we construct a Hochschild cohomology complex for $A$ and use this to determine the Hochschild cohomology of $A$ in the cases listed above.

The Hochschild (co)homology and cyclic homology of $A$ can be regarded as invariants for the braided subfactors associated to the $SU(3)$ modular invariants.
Beginning with a pair $(\mathcal{G},W)$ given by a cell system $W$ on an $SU(3)$ $\mathcal{ADE}$ graph $\mathcal{G}$, we construct a braided subfactor $N \subset M$ which yields a nimrep which recovers the graph $\mathcal{G}$ as described above. Then we can construct the algebra $A(\mathcal{G},W)$ whose Hochschild (co)homology and cyclic homology only depends on the original pair $(\mathcal{G},W)$, or equivalently, on the braided subfactor $N \subset M$.

\section{Almost Calabi-Yau algebras} \label{sect:almostCYalg}

Let $\mathcal{G}$ be a finite directed graph, and denote by $\mathcal{G}_n$ the set of all paths on $\mathcal{G}$ of length $n$. The vertices of $\mathcal{G}$ are the paths of length 0.
If $a \in \mathcal{G}_1$ is an edge on $\mathcal{G}$, we denote by $\widetilde{a} \in \mathcal{G}^{\mathrm{op}}_1$ the corresponding edge with opposite orientation on $\mathcal{G}^{\mathrm{op}}$.
The path algebra $\mathbb{C}\mathcal{G} = \bigoplus_{k=0}^{\infty} (\mathbb{C}\mathcal{G})_k$ is the graded complex vector space with basis of the $k^{\mathrm{th}}$-graded part $(\mathbb{C}\mathcal{G})_k$ given by $\mathcal{G}_k$, where paths may begin at any vertex of $\mathcal{G}$. Multiplication of two paths $a \in (\mathbb{C}\mathcal{G})_k$ and $b \in (\mathbb{C}\mathcal{G})_l$ is given by concatenation of paths $a \cdot b \in (\mathbb{C}\mathcal{G})_{k+l}$ (or simply $ab$), with $ab$ defined to be zero if $r(a) \neq s(b)$, where $s(a)$, $r(a)$ denotes the source, range vertex respectively of the path $a$.
The commutator quotient $\mathbb{C}\mathcal{G} / [\mathbb{C}\mathcal{G}, \mathbb{C}\mathcal{G}]$ may be identified, up to cyclic permutation of the arrows, with the vector space spanned by cyclic paths in $\mathcal{G}$.
Let $\partial_a : \mathbb{C}\mathcal{G} / [\mathbb{C}\mathcal{G}, \mathbb{C}\mathcal{G}] \rightarrow \mathbb{C}\mathcal{G}$ be the derivation given by
$\partial_{a} (a_1 \cdots a_n) = \sum_{j} a_{j+1} \cdots a_n a_1 \cdots a_{j-1}$,
where the summation is over all indices $j$ such that $a_j = a$.
Then for a potential $\Phi \in \mathbb{C}\mathcal{G} / [\mathbb{C}\mathcal{G}, \mathbb{C}\mathcal{G}]$, which is some linear combination of cyclic paths in $\mathcal{G}$, we define the algebra
$A(\mathbb{C}\mathcal{G}, \Phi) = \mathbb{C}\mathcal{G} / \{ \partial_a \Phi \}$,
which is the quotient of the path algebra by the two-sided ideal generated by the relations $\partial_a \Phi \in \mathbb{C}\mathcal{G}$, for all edges $a$ of $\mathcal{G}$.
The Hilbert series $H_A$ for $A(\mathbb{C}\mathcal{G}, \Phi)$ is defined as $H_A(t) = \sum_{k=0}^{\infty} H_{ji}^k t^k$, where the $H_{ji}^k$ are matrices which count the dimension of the subspace $\{ i x j | \; x \in A(\mathbb{C}\mathcal{G}, \Phi)_k \}$, where $A(\mathbb{C}\mathcal{G}, \Phi)_k$ is the subspace of $A(\mathbb{C}\mathcal{G}, \Phi)$ of all paths of length $k$, and $i,j \in A(\mathbb{C}\mathcal{G}, \Phi)_0$.

\begin{figure}[tb]
\begin{center}
\includegraphics[width=100mm]{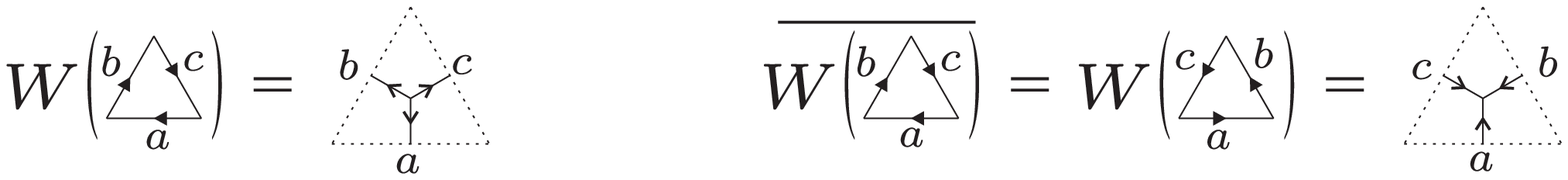}\\
 \caption{Cells associated to trivalent vertices} \label{fig:Oc-Kup}
\end{center}
\end{figure}

Ocneanu \cite{ocneanu:2000ii} defined a cell system $W$ on any $SU(3)$ $\mathcal{ADE}$ graph $\mathcal{G}$, associating a complex number $W \left( \triangle_{i,j,k}^{(a,b,c)} \right)$, now called an Ocneanu cell, to each closed loop of length three $\triangle_{i,j,k}^{(a,b,c)}$ in $\mathcal{G}$ as in Figure \ref{fig:Oc-Kup}, where $a,b,c$ are edges on $\mathcal{G}$, and $i,j,k$ are the vertices on $\mathcal{G}$ given by $i=s(a)=r(c)$, $j=s(b)=r(a)$, $k=s(c)=r(b)$.
These cells satisfy two properties, called Ocneanu's type I, II equations respectively, which are obtained by evaluating the Kuperberg relations K2, K3 for an $A_2$-spider \cite{kuperberg:1996} using the identification in Figure \ref{fig:Oc-Kup}: \\
$(i)$ for any type I frame \includegraphics[width=16mm]{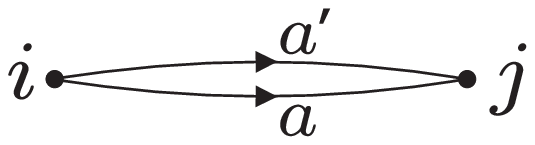} in $\mathcal{G}$ we have
$$\sum_{k,b_1,b_2} W \left( \triangle_{i,j,k}^{(a,b_1,b_2)} \right) \overline{W \left( \triangle_{i,j,k}^{(a',b_1,b_2)} \right)} = \delta_{a,a'} [2]_q \phi_i \phi_j$$
$(ii)$ for any type II frame \includegraphics[width=30mm]{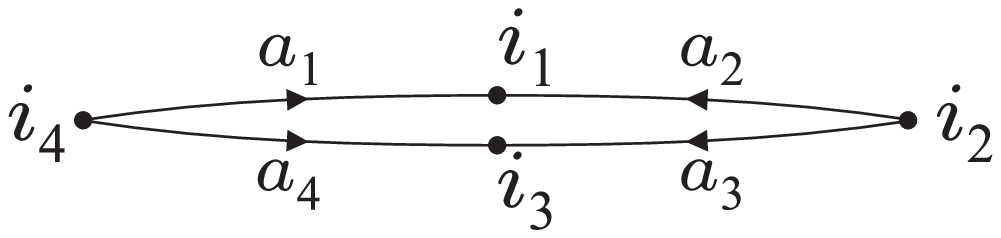} in $\mathcal{G}$ we have
\begin{eqnarray*}
\lefteqn{ \sum_{k,b_j} \phi_k^{-1} W \left( \triangle_{i_2,i_1,k}^{(a_2,b_1,b_2)} \right) \overline{W \left( \triangle_{i_2,i_3,k}^{(a_3,b_3,b_2)} \right)} W \left( \triangle_{i_4,i_3,k}^{(a_4,b_3,b_4)} \right) \overline{W \left( \triangle_{i_4,i_1,k}^{(a_1,b_4,b_1)} \right)} } \nonumber \\
& \qquad & = \delta_{a_1,a_4} \delta_{a_2,a_3} \phi_{i_4} \phi_{i_1} \phi_{i_2} + \delta_{a_1,a_2} \delta_{a_3,a_4} \phi_{i_1} \phi_{i_2} \phi_{i_3} \hspace{25mm} \label{eqn:typeII_frame}
\end{eqnarray*}
Here $(\phi_v)_v$ is the Perron-Frobenius eigenvector for the Perron-Frobenius eigenvalue $\alpha = [3]_q$ of $\mathcal{G}$.
The existence of these cells for the finite $\mathcal{ADE}$ graphs was claimed by Ocneanu \cite{ocneanu:2000ii}, and shown in \cite{evans/pugh:2009i} with the exception of the graph $\mathcal{E}_4^{(12)}$. These cells define a unitary connection on the graph $\mathcal{G}$ which satisfy the Yang-Baxter equation \cite[Lemma 3.2]{evans/pugh:2009i}.

Two cell systems $W_1$, $W_2$ on an $SU(3)$ $\mathcal{ADE}$ graph $\mathcal{G}$ are equivalent if, for each pair of adjacent vertices $i$, $j$ of $\mathcal{G}$, we can find a family of unitary matrices $(u(a,b))_{a,b}$, where $a$, $b$ are any pair of edges from $i$ to $j$, such that
$$W_1(\triangle_{i_1,i_2,i_3}^{(a_1,a_2,a_3)}) = \sum_{a_1',a_2',a_3'} u(a_1,a_1') u(a_2,a_2') u(a_3,a_3') W_2(\triangle_{i_1,i_2,i_3}^{(a_1',a_2',a_3')}),$$
where $a_l$ are edges from $i_l$ to $i_{l+1}$, and the sum is over all edges $a_l'$ from $i_l$ to $i_{l+1}$, $l=1,2,3$.

There is up to equivalence precisely one connection on the graphs $\mathcal{A}^{(m)}$, $\mathcal{A}^{(2m+1)\ast}$, $\mathcal{E}^{(8)}$, $\mathcal{E}^{(8)\ast}$, $\mathcal{E}_5^{(12)}$ and $\mathcal{E}^{(24)}$. For the graphs $\mathcal{A}^{(2m)\ast}$ and $\mathcal{E}_2^{(12)}$ there are precisely two inequivalent connections, which are obtained from each other by a $\mathbb{Z}_2$ symmetry of the graph. This $\mathbb{Z}_2$ symmetry is the conjugation of the graph in the case of $\mathcal{E}_2^{(12)}$. There is at least one connection for each graph $\mathcal{D}^{(m)}$, $m \not \equiv 0 \textrm{ mod } 3$, and at least two inequivalent connections for each graph $\mathcal{D}^{(3p)}$, which are the complex conjugates of each other. There is at least one connection for each graph $\mathcal{D}^{(2m+1)\ast}$, and at least two inequivalent connections for each graph $\mathcal{D}^{(2m)\ast}$, which are obtained from each other by a $\mathbb{Z}_2$ symmetry of the graph. There are also at least two inequivalent connections for the graph $\mathcal{E}_1^{(12)}$, which are obtained from each other by conjugation of the graph.

For the $SU(3)$ $\mathcal{ADE}$ graphs, we define the almost Calabi-Yau algebra $A(\mathcal{G},W)$ to be the graded quotient algebra
$$A(\mathcal{G},W) := A(\mathbb{C} \mathcal{G}, \Phi_W),$$
where the potential $\Phi_W$ is given by \cite[equation (40)]{evans/pugh:2010ii} (see also \cite[Remark 4.5.7]{ginzburg:2006}):
$$\Phi_W = \sum_{abc} W(\triangle_{abc}) \triangle_{abc} \quad \in \mathbb{C} \mathcal{G} / [\mathbb{C} \mathcal{G}, \mathbb{C} \mathcal{G}],$$
where the summation is over all closed paths $abc$ of length 3 on $\mathcal{G}$. The grading on $\mathbb{C}\mathcal{G}$ descends to the quotient algebra $A = A(\mathcal{G},W)$.
These almost Calabi-Yau algebras were studied in \cite{evans/pugh:2010ii} for all the cell systems constructed in \cite{evans/pugh:2009i}.
Equivalent cell systems yield isomorphic almost Calabi-Yau algebras.
For any cell system $W$, we can take its complex conjugate $\overline{W}$ to obtain another (possibly equivalent) cell system. The almost Calabi-Yau algebra for $\overline{W}$ is isomorphic to that for $W$.
The conjugation $\tau: {}_N \mathcal{X}_N \rightarrow {}_N \mathcal{X}_N$ on the braided system of endomorphisms of $SU(3)_k$ on a factor $N$, given by the conjugation on the representations of $SU(3)$, induces a conjugation $\tau: {}_N \mathcal{X}_M \rightarrow {}_N \mathcal{X}_M$ such that $G_{\overline{\lambda}} = \tau G_{\lambda} \tau$, where $G_{\lambda} a = \lambda a$ for $\lambda \in {}_N \mathcal{X}_N$, $a \in {}_N \mathcal{X}_M$. For any cell system $W = W^+$, this conjugation of the graph yields a conjugate cell system $W^-$, which might be equivalent to $W^+$. The almost Calabi-Yau algebra for $W^-$ is anti-isomorphic to that for $W^+$.

The Hilbert series $H_A(t)$ of $A(\mathcal{G},W)$, for an $SU(3)$ $\mathcal{ADE}$ graph $\mathcal{G}$ with adjacency matrix $\Delta_{\mathcal{G}}$, Coxeter number $h=k+3$ and cell system $W$, is given by \cite[Theorem 3.1]{evans/pugh:2010ii}
\begin{equation} \label{eqn:Hilbert_Series-SU(3)ADE}
H_A (t) = \frac{1 - P t^h}{1 - \Delta_{\mathcal{G}} t + \Delta_{\mathcal{G}}^T t^2 - t^3},
\end{equation}
where $P$ is the permutation matrix corresponding to a $\mathbb{Z}_3$ symmetry of the graph.
It is the identity for $\mathcal{D}^{(n)}$, $\mathcal{A}^{(n)\ast}$, $n \geq 5$, $\mathcal{E}^{(8)\ast}$, $\mathcal{E}_l^{(12)}$, $l=1,2,4,5$, and $\mathcal{E}^{(24)}$. For the remaining graphs $\mathcal{A}^{(n)}$, $\mathcal{D}^{(n) \ast}$ and $\mathcal{E}^{(8)}$, let $V$ be the permutation matrix corresponding to the clockwise rotation of the graph by $2 \pi /3$. Then
$$ P = \left\{
\begin{array}{cl} V^2 & \mbox{ for } \quad \mathcal{A}^{(n)}, n \geq 4, \\
                  V & \mbox{ for } \quad \mathcal{E}^{(8)}, \\
                  V^{2n} & \mbox{ for } \quad \mathcal{D}^{(n) \ast}, n \geq 5.
\end{array} \right.$$
The numerator and denominator in (\ref{eqn:Hilbert_Series-SU(3)ADE}) commute, since any permutation matrix which corresponds to a symmetry of the graph $\mathcal{G}$ commutes with $\Delta_{\mathcal{G}}$ and $\Delta_{\mathcal{G}}^T$.

\subsection{Periodic resolution for almost Calabi-Yau algebras} \label{sect:resolution_almostCYalg}

We define a non-degenerate form on $A$ by setting $f$ to be the function which is 0 on every element of $A$ of length $< h-3$, and 1 on $u_{i\nu(i)}$ for some $i \in \mathcal{G}_1$, where $u_{j\nu(j)}$ denotes a generator of the one-dimensional top-degree space $j \cdot A_{h-3} \cdot \nu(j)$, where $\nu$ is the permutation of the vertices of $\mathcal{G}$ given by the permutation matrix $P$ in (\ref{eqn:Hilbert_Series-SU(3)ADE}). Then using the relation $(x,y) = (y,\beta(x))$ this determines the value of $f$ on $u_{j\nu(j)}$, for all other $j \in \mathcal{G}_1$. We normalize the $u_{j\nu(j)}$ such that $f(u_{j\nu(j)}) = 1$ for all $j \in \mathcal{G}_1$.
The image of the simple object $f_{(k,0)} \in A_2\textrm{-}TL^{(k)}$ under the functor $F$ given by (\ref{eqn:functorF}) defines a unique permutation $\nu$ of the graph $\mathcal{G}$, which is described as follows.
The permutation $\nu$ of the graph is given by the $\mathbb{Z}_3$ symmetry which defines the permutation matrix $P$ in (\ref{eqn:Hilbert_Series-SU(3)ADE}) (note that there are no double edges on the graphs $\mathcal{G}$ for which $P$ is non-trivial).
Then the Nakayama automorphism $\beta$ of $A$ is defined on $\mathcal{G}$ by $\beta = \nu$ \cite[Theorem 4.6]{evans/pugh:2010ii}.

Now $A$ has the following finite resolution as an $A$-$A$ bimodule \cite[Theorem 5.1]{evans/pugh:2010ii}:
\begin{equation} \label{exact_seq-almostCY}
0 \rightarrow \mathcal{N}[h] \stackrel{\iota_0}{\rightarrow} A \otimes_S A[3] \stackrel{\mu_3}{\rightarrow} A \otimes_S \widetilde{V} \otimes_S A[1] \stackrel{\mu_2}{\rightarrow} A \otimes_S V \otimes_S A \stackrel{\mu_1}{\rightarrow} A \otimes_S A \stackrel{\mu_0}{\rightarrow} A \rightarrow 0.
\end{equation}
Here $S$ is the $A$-$A$ bimodule $(\mathbb{C}\mathcal{G})_0$, and $V$, $\widetilde{V}$ are the $A$-$A$ bimodules generated by $\mathcal{G}_1$, $\mathcal{G}^{\mathrm{op}}_1$ respectively. The $A$-$A$ bimodule $\mathcal{N} = {}_1 A_{\beta^{-1}}$ is equal to $A$ as a vector space. The left $A$-action is given by concatenation, but the right $A$-action is twisted by the inverse of the Nakayama automorphism $\beta$, i.e. $a \cdot x \cdot b = ax\beta^{-1}(b)$ for all $a,b \in A$, $x \in \mathcal{N}$.
The connecting $A$-$A$ bimodule maps are given by
\begin{eqnarray}
\mu_0(1 \otimes 1) & = & 1, \label{mu_0} \\
\mu_1(1 \otimes a \otimes 1) & = & a \otimes 1 - 1 \otimes a, \label{mu_1} \\
\mu_2(1 \otimes \widetilde{a} \otimes 1) & = & \sum_{b,b' \in \mathcal{G}_1} W_{abb'} (b \otimes b' \otimes 1 + 1 \otimes b \otimes b'), \label{mu_2} \\
\mu_3(1 \otimes 1) & = & \sum_{a \in \mathcal{G}_1} a \otimes \widetilde{a} \otimes 1 - \sum_{a \in \mathcal{G}_1} 1 \otimes \widetilde{a} \otimes a, \label{mu_3} \\
\iota_0(1) & = & \sum_j w_j \otimes w_j^{\ast}, \nonumber
\end{eqnarray}
where $\{ w_j \}$ is a homogeneous basis for $A$, and $\{ w_j^{\ast} \}$ is its corresponding dual basis, i.e. $w_j w_j^{\ast} = u_{i \nu(i)}$ where $i = s(w_j)$. The $A$-$A$ bimodule $B = B^{(1)} \otimes_S \cdots \otimes_S B^{(p)}$ is equipped with the \emph{total grading} which comes from the grading on the graded $A$-$A$ bimodules $B^{(i)}$, that is, $B = \bigoplus_{k=0}^{\infty} B_k$ where $B_k = \bigoplus_{k_i: \sum_{i=1}^p k_i = k} B^{(1)}_{k_1} \otimes_S \cdots \otimes_S B^{(p)}_{k_p}$.

For each $SU(3)$ $\mathcal{ADE}$ graph, the Nakayama automorphism has order 3, $\beta^{3} = \mathrm{id}$, so we can make a canonical identification $A = \mathcal{N} \otimes_{A} \mathcal{N} \otimes_{A} \mathcal{N}$. We let $\mathcal{N}^{(k)} := {}_1 A_{\beta^{-k}}$, for $k \in \mathbb{Z}$. In particular, we have $A = \mathcal{N}^{(0)}$, $\mathcal{N} = \mathcal{N}^{(1)}$ and $\mathcal{N}^{(2)} = {}_1 A_{\beta} = \mathcal{N} \otimes_{A} \mathcal{N}$. Note that for graphs with trivial Nakayama automorphism, $A = \mathcal{N}^{(k)}$ as $A$-$A$ bimodules, for all $k \in \mathbb{Z}$.

Applying the functor $- \otimes_{A} \mathcal{N}$ to the exact sequence (\ref{exact_seq-almostCY}) we obtain the exact sequence:
\begin{eqnarray*}
& 0 \rightarrow \mathcal{N}^{(2)}[2h] \stackrel{\iota_2}{\rightarrow} A \otimes_S \mathcal{N}[h+3] \stackrel{\mu_7}{\rightarrow} A \otimes_S \widetilde{V} \otimes_S \mathcal{N}[h+1] \stackrel{\mu_6}{\rightarrow} A \otimes_S V \otimes_S \mathcal{N}[h] & \nonumber \\
& \stackrel{\mu_5}{\rightarrow} A \otimes_S \mathcal{N}[h] \stackrel{\iota_1}{\rightarrow} \mathcal{N}[h] \rightarrow 0, & \label{exact_seq-almostCY-2}
\end{eqnarray*}
where $\iota_1(x \otimes y) = xy$, $\iota_2(a) = a \sum_j w_j \otimes w_j^{\ast}$, where $\{ w_j \}$ is a homogeneous basis for $A$ and $\{ w_j^{\ast} \}$ is its corresponding dual basis, and $\mu_{i+4} = \mu_i$. Similarly, applying the functor a second time we obtain the exact sequence:
\begin{eqnarray*}
& 0 \rightarrow A[3h] \stackrel{\iota_4}{\rightarrow} A \otimes_S \mathcal{N}^{(2)}[2h+3] \stackrel{\mu_{11}}{\rightarrow} A \otimes_S \widetilde{V} \otimes_S \mathcal{N}^{(2)}[2h+1] \stackrel{\mu_{10}}{\rightarrow} A \otimes_S V \otimes_S \mathcal{N}^{(2)}[2h] & \nonumber \\
& \stackrel{\mu_9}{\rightarrow} A \otimes_S \mathcal{N}^{(2)}[2h] \stackrel{\iota_3}{\rightarrow} \mathcal{N}^{(2)}[2h] \rightarrow 0. & \label{exact_seq-almostCY-3}
\end{eqnarray*}

We now construct a projective resolution of $A$, that is, a resolution of $A$ by projective modules.
Setting $\mu_4 = \iota_0 \iota_1$, $\mu_8 = \iota_2 \iota_3$, we obtain the following projective resolution of $A$, which is periodic with period 12:
\begin{eqnarray}
& \cdots \,\, \rightarrow A[3h] \stackrel{\mu_{12}}{\rightarrow} A \otimes_S \mathcal{N}^{(2)}[2h+3] \stackrel{\mu_{11}}{\rightarrow} A \otimes_S \widetilde{V} \otimes_S \mathcal{N}^{(2)}[2h+1] \stackrel{\mu_{10}}{\rightarrow} A \otimes_S V \otimes_S \mathcal{N}^{(2)}[2h] & \nonumber \\
& \stackrel{\mu_9}{\rightarrow} A \otimes_S \mathcal{N}^{(2)}[2h] \stackrel{\mu_8}{\rightarrow} A \otimes_S \mathcal{N}[h+3] \stackrel{\mu_7}{\rightarrow} A \otimes_S \widetilde{V} \otimes_S \mathcal{N}[h+1] \stackrel{\mu_6}{\rightarrow} A \otimes_S V \otimes_S \mathcal{N}[h] & \nonumber \\
& \stackrel{\mu_5}{\rightarrow} A \otimes_S \mathcal{N}[h] \stackrel{\mu_4}{\rightarrow} A \otimes_S A[3] \stackrel{\mu_3}{\rightarrow} A \otimes_S \widetilde{V} \otimes_S A[1] \stackrel{\mu_2}{\rightarrow} A \otimes_S V \otimes_S A \stackrel{\mu_1}{\rightarrow} A \otimes_S A \stackrel{\mu_0}{\rightarrow} A \rightarrow 0, & \nonumber \\
&& \label{resolution-almostCY}
\end{eqnarray}
where the connecting maps $\mu_i$ are given by (\ref{mu_0})-(\ref{mu_3}) for $0 \leq i \leq 3$, $\mu_4(x \otimes y) = xy \sum_j w_j \otimes w_j^{\ast}$, where $\{ w_j \}$ is a homogeneous basis for $A$ and $\{ w_j^{\ast} \}$ is its corresponding dual basis, and $\mu_i = \mu_{i-4}$ for $i \geq 5$.

Thus we find that the Hochschild (co)homology of $A$ is periodic with period 12, i.e. the grading is shifted by $3h$ ($-3h$) when the degree of the homology (respectively cohomology) is shifted by 12.
In the case of trivial Nakayama automorphism the Hochschild (co)homology of $A$ in fact has period 4.

\section{The Hochschild homology of $A(\mathcal{G},W)$} \label{sect:Hoch_hom}
\subsection{The Hochschild homology complex} \label{sect:Hoch_hom-complex}

In this section we will construct a complex which determines the Hochschild homology of the almost Calabi-Yau algebra $A = A(\mathcal{G},W)$.

Let $A^{\mathrm{op}}$ denote the algebra with opposite multiplication, i.e. $a \cdot b = ba$, and define $A^e = A^{\mathrm{op}} \otimes_S A$. Any $A$-$A$ bimodule becomes a left $A^e$-module, and vice versa, by defining the left action of $A^e$ on $A$ by $(a \otimes b) x = bxa$ for all $x \in A$, $a \otimes b \in A^{\mathrm{op}} \otimes_S A$.

The Hochschild homology $HH_{\bullet}(A)$ of $A$ may be defined to be the derived functor $HH_n(A) = \mathrm{Tor}_n^{A^e}(A,A)$, e.g. \cite[Proposition 1.1.13]{loday:1998}, i.e. as the homology of the complex
$$\cdots \rightarrow P_2 \otimes_{A^e} A \rightarrow P_1 \otimes_{A^e} A \rightarrow P_0 \otimes_{A^e} A \rightarrow A \otimes_{A^e} A \rightarrow 0$$
where $\quad \cdots \rightarrow P_2 \rightarrow P_1 \rightarrow P_0 \rightarrow A \rightarrow 0$ is any projective resolution of $A$.

For an $A$-$A$ bimodule $M$, denote by $M^S$ the $S$-centralizer sub-bimodule given by all elements $x \in M$ such that $i x = x i$ for all $i \in S$.
We make the following identifications, for $k=0,1,2$ (c.f. \cite{etingof/eu:2007}):
$$\begin{array}{ll}
(A \otimes_S \mathcal{N}^{(k)}) \otimes_{A^e} A = (\mathcal{N}^{(k)})^S: & (x \otimes y) \otimes z = y\beta^{-k}(zx), \\
(A \otimes_S V \otimes_S \mathcal{N}^{(k)}) \otimes_{A^e} A = (V \otimes_S \mathcal{N}^{(k)})^S: & (x \otimes a \otimes y) \otimes z = a \otimes y\beta^{-k}(zx), \\
(A \otimes_S \widetilde{V} \otimes_S \mathcal{N}^{(k)}) \otimes_{A^e} A = (\widetilde{V} \otimes_S \mathcal{N}^{(k)})^S: & (x \otimes \widetilde{a} \otimes y) \otimes z = \widetilde{a} \otimes y\beta^{-k}(zx),
\end{array}$$
where the left and right hand sides have the same total degree.
Thus, applying the functor $- \otimes_{A^e} A$ to the resolution (\ref{resolution-almostCY}), we obtain the Hochschild homology complex:
\begin{eqnarray}
& \cdots \,\, \rightarrow A^S[3h] \stackrel{\mu_{12}'}{\rightarrow} (\mathcal{N}^{(2)})^S[2h+3] \stackrel{\mu_{11}'}{\rightarrow} (\widetilde{V} \otimes_S \mathcal{N}^{(2)})^S[2h+1] \stackrel{\mu_{10}'}{\rightarrow} (V \otimes_S \mathcal{N}^{(2)})^S[2h] & \nonumber \\
& \stackrel{\mu_9'}{\rightarrow} (\mathcal{N}^{(2)})^S[2h] \stackrel{\mu_8'}{\rightarrow} \mathcal{N}^S[h+3] \stackrel{\mu_7'}{\rightarrow} (\widetilde{V} \otimes_S \mathcal{N})^S[h+1] \stackrel{\mu_6'}{\rightarrow} (V \otimes_S \mathcal{N})^S[h] & \nonumber \\
& \stackrel{\mu_5'}{\rightarrow} \mathcal{N}^S[h] \stackrel{\mu_4'}{\rightarrow} A^S[3] \stackrel{\mu_3'}{\rightarrow} (\widetilde{V} \otimes_S A)^S[1] \stackrel{\mu_2'}{\rightarrow} (V \otimes_S A)^S \stackrel{\mu_1'}{\rightarrow} A^S \rightarrow 0, & \label{HH_hom_complex-almostCY}
\end{eqnarray}
where the connecting maps are given, for $k=0,1,2,\ldots \;\;$ by
\begin{eqnarray*}
\mu_{4k+1}'(a \otimes x) & = & \mu_{4k+1}(1 \otimes a \otimes 1) \otimes_{A^e} \beta^{k}(x) = (a \otimes 1 - 1 \otimes a) \otimes_{A^e} \beta^{k}(x) \\
& = & x \beta^{-k}(a) - ax, \\
\mu_{4k+2}'(\widetilde{a} \otimes x) & = & \mu_{4k+2}(1 \otimes \widetilde{a} \otimes 1) \otimes_{A^e} \beta^{k}(x) \\
& = & \sum_{b,b' \in \mathcal{G}_1} W_{abb'} (b \otimes b' \otimes 1 + 1 \otimes b \otimes b') \otimes_{A^e} \beta^{k}(x) \\
& = & \sum_{b,b' \in \mathcal{G}_1} W_{abb'} (b' \otimes x \beta^{-k}(b) + b \otimes b'x), \\
\mu_{4k+3}'(x) & = & \mu_{4k+3}(1 \otimes 1) \otimes_{A^e} \beta^{k}(x) = \left( \sum_{a \in \mathcal{G}_1} a \otimes \widetilde{a} \otimes 1 - 1 \otimes \widetilde{a} \otimes a \right) \otimes_{A^e} \beta^{k}(x) \\
& = & \sum_{a \in \mathcal{G}_1} \widetilde{a} \otimes (x \beta^{-k}(a) - ax), \\
\mu_{4k+4}'(y) & = & \mu_{4k+4}(1 \otimes 1) \otimes_{A^e} \beta^{k+1}(y) = \left( \sum_j w_j \otimes w_j^{\ast} \right) \otimes_{A^e} \beta^{k+1}(y) \\
& = & \sum_j w_j^{\ast} \beta(y) \beta^{-k}(w_j),
\end{eqnarray*}
where $a \in V$, $x \in \mathcal{N}^{(k)}$, $y \in \mathcal{N}^{(k+1)}$, $\{ w_j \}$ is a homogeneous basis for $A$ and $\{ w_j^{\ast} \}$ is its corresponding dual basis.

We will now show that this complex has a self-duality.
Using the non-degenerate form, we can make the identifications $\mathcal{N}^{(k)} = (\mathcal{N}^{(2-k)})^{\ast}[h-3]$ by sending $x \mapsto (-,x)$.
We can define a non-degenerate form on $(V \oplus \widetilde{V}) \otimes_S \mathcal{N}^{(k)}$ by $(a_1 \otimes x_1, a_2 \otimes x_2) = \delta_{a_1,\beta^{k-1}(\widetilde{a_2})} (x_1,x_2)$ for $x_1 \in \mathcal{N}^{(2-k)}$, $x_2 \in \mathcal{N}^{(k)}$, and $a_1 \in V_1$, $a_2 \in V_2$, where $V_i \in \{ V, \widetilde{V}\}$, $i=1,2$. For the $\mathcal{A}^{\ast}$ graphs, $V = \widetilde{V}$ and we replace $(V \oplus \widetilde{V}) \otimes_S \mathcal{N}^{(k)}$ above by $V \otimes_S \mathcal{N}^{(k)}$.
This allows us to make identifications $V \otimes_S \mathcal{N}^{(k)} = (\widetilde{V} \otimes_S \mathcal{N}^{(2-k)})^{\ast}[h-1]$, $\widetilde{V} \otimes_S \mathcal{N}^{(k)} = (V \otimes_S \mathcal{N}^{(2-k)})^{\ast}[h-1]$, by sending $a \otimes x \mapsto (-,a \otimes x)$.

If we take the Hochschild homology sequence (\ref{HH_hom_complex-almostCY}) and dualise, we get:
\begin{eqnarray*}
& \cdots \,\, \stackrel{(\mu_{12}')^{\ast}}{\leftarrow} A^S[-3h] \stackrel{(\mu_{11}')^{\ast}}{\leftarrow} (V \otimes_S A)^S[-3h] \stackrel{(\mu_{10}')^{\ast}}{\leftarrow} (\widetilde{V} \otimes_S A)^S[-3h+1] \stackrel{(\mu_9')^{\ast}}{\leftarrow} & \\
& \stackrel{(\mu_9')^{\ast}}{\leftarrow} A^S[-3h+3] \stackrel{(\mu_8')^{\ast}}{\leftarrow} \mathcal{N}^S[-2h] \stackrel{(\mu_7')^{\ast}}{\leftarrow} (V \otimes_S \mathcal{N})^S[-2h] \stackrel{(\mu_6')^{\ast}}{\leftarrow} & \\
& \stackrel{(\mu_6')^{\ast}}{\leftarrow} (\widetilde{V} \otimes_S \mathcal{N})^S[-2h+1] \stackrel{(\mu_5')^{\ast}}{\leftarrow} \mathcal{N}^S[-2h+3] \stackrel{(\mu_4')^{\ast}}{\leftarrow} (\mathcal{N}^{(2)})^S[-h] \stackrel{(\mu_3')^{\ast}}{\leftarrow} & \\
& \stackrel{(\mu_3')^{\ast}}{\leftarrow} (V \otimes_S \mathcal{N}^{(2)})^S[-h] \stackrel{(\mu_2')^{\ast}}{\leftarrow} (\widetilde{V} \otimes_S \mathcal{N}^{(2)})^S[-h+1] \stackrel{(\mu_1')^{\ast}}{\leftarrow} (\mathcal{N}^{(2)})^S[-h+3] \leftarrow 0. &
\end{eqnarray*}

\begin{Prop}
We have $\mu_i' = \pm(\mu_{12-i}')^{\ast}$, $i=1,\ldots,11$.
\end{Prop}

\noindent \emph{Proof}:
(i) $\mu_1' = - (\mu_{11}')^{\ast}$:
Let $a \in V$, $x \in A$ and $y \in \mathcal{N}^{(2)}$. Then
\begin{eqnarray*}
(\mu_1'(a \otimes x),y) & = & (xa-ax,y) = (x,ay-y\beta(a)) = (a \otimes x, -\sum_{b \in \mathcal{G}_1} \widetilde{b} \otimes (y\beta(b)-by)) \\
& = & (a \otimes x,-\mu_{11}'(y)).
\end{eqnarray*}
(ii) $\mu_2' = (\mu_{10}')^{\ast}$:
Let $a,a' \in V$, $x \in A$ and $y \in \mathcal{N}^{(2)}$. Then
\begin{eqnarray*}
\lefteqn{ \hspace{-10mm} (\mu_2'(\widetilde{a} \otimes x),\widetilde{a'} \otimes y) \;\; = \;\; (\sum_{b,b' \in \mathcal{G}_1} W_{abb'} (b' \otimes xb + b \otimes b'x), \widetilde{a'} \otimes y) } \\
\hspace{10mm} & = & (\sum_{b \in \mathcal{G}_1} W_{aba'} xb + \sum_{b' \in \mathcal{G}_1} W_{aa'b'} b'x,y) \;\; = \;\; (x, \sum_{b \in \mathcal{G}_1} W_{aba'} by + \sum_{b' \in \mathcal{G}_1} W_{aa'b'} y\beta(b')) \\
& = & (\widetilde{a} \otimes x, \sum_{b,b' \in \mathcal{G}_1} W_{b'ba'} (b' \otimes by + b \otimes y\beta(b')) = (\widetilde{a} \otimes x, \mu_{10}'(\widetilde{a'} \otimes y)).
\end{eqnarray*}
(iii) $\mu_3' = - (\mu_9')^{\ast}$:
Let $a' \in V$, $x \in A$ and $y \in \mathcal{N}^{(2)}$. Then
\begin{eqnarray*}
(\mu_3'(x),a' \otimes y) & = & (\sum_{a \in \mathcal{G}_1} \widetilde{a} \otimes (xa-ax),a' \otimes y) = (xa'-a'x,y) = (x,a'y-y\beta(a')) \\
& = & (x,-\mu_9'(a' \otimes y)).
\end{eqnarray*}
(iv) $\mu_4' = (\mu_8')^{\ast}$:
Let $x \in \mathcal{N}$ and $y \in \mathcal{N}^{(2)}$. Then
\begin{eqnarray*}
(\mu_4'(x),y) & = & (\sum_j w_j^{\ast}\beta(x)w_j,y) = (\sum_j w_j x \beta^2(w_j^{\ast}),y) = (x,\sum_j \beta^2(w_j^{\ast})y\beta(w_j)) \\
& = & (x,\sum_j w_j^{\ast} \beta(y) \beta^2(w_j)) = (x,\mu_8'(y)),
\end{eqnarray*}
where the second equality holds since if $\{ w_j^{\ast} \}$ is a dual basis of $\{ w_j \}$, then $\{ w_j \}$ is a dual basis of $\{ \beta^2(w_j^{\ast}) \}$, and $\sum_j w_j^{\ast}\beta(x)w_j = 0 = \sum_j w_j x \beta^2(w_j^{\ast})$ unless $|x|=0$ such that $\beta(x)=x$. The penultimate equality is given by replacing the basis $\{ \beta^2(w_j^{\ast}) \}$ with the equivalent basis $\{ w_j^{\ast} \}$, and the fact that $\sum_j \beta^2(w_j^{\ast})y\beta(w_j) = 0 = \sum_j w_j^{\ast} \beta(y) \beta^2(w_j)$ unless $|y|=0$ such that $\beta(y)=y$. \\
(v) $\mu_5' = - (\mu_7')^{\ast}$:
Let $a \in V$ and $x, y \in \mathcal{N}$. Then
\begin{eqnarray*}
(\mu_5'(a \otimes x),y) & = & (x\beta^2(a)-ax,y) = (x,\beta^2(a)y-y\beta(a)) \\
& = & (a \otimes x, -\sum_{b \in \mathcal{G}_1} \widetilde{b} \otimes (y\beta^2(b)-by)) = (a \otimes x,-\mu_7'(y)).
\end{eqnarray*}
(vi) $\mu_6' = (\mu_6')^{\ast}$:
Let $a,a' \in V$ and $x, y \in \mathcal{N}$. Then
\begin{eqnarray*}
\lefteqn{ \hspace{-2mm} (\mu_6'(\widetilde{a} \otimes x),\widetilde{a'} \otimes y) \;\; = \;\; (\sum_{b,b' \in \mathcal{G}_1} W_{abb'} (b' \otimes x\beta^2(b) + b \otimes b'x), \widetilde{a'} \otimes y) } \\
& = & (\sum_{b \in \mathcal{G}_1} W_{aba'} x\beta^2(b) + \sum_{b' \in \mathcal{G}_1} W_{aa'b'} b'x,y) \;\; = \;\; (x, \sum_{b \in \mathcal{G}_1} W_{aba'} \beta^2(b)y + \sum_{b' \in \mathcal{G}_1} W_{aa'b'} y\beta(b')) \\
& = & (\widetilde{a} \otimes x, \sum_{b,b' \in \mathcal{G}_1} W_{b'ba'} (b' \otimes by + b \otimes y\beta^2(b')) = (\widetilde{a} \otimes x, \mu_6'(\widetilde{a'} \otimes y)).
\end{eqnarray*}

\vspace{-5mm} $\,$
\hfill
$\Box$

Note however that $(\mu_{12}')^{\ast} = \mu_{12}' \circ \beta$:
Let $x,y \in A$. Then
\begin{eqnarray*}
(\mu_{12}'(x),y) & = & (\sum_j w_j^{\ast}\beta(x)\beta(w_j),y) = (\beta(x),\sum_j \beta(w_j) y \beta(w_j^{\ast})) = (x,\sum_j w_j\beta^2(y)w_j^{\ast}) \\
& = & (x,\sum_j w_j^{\ast}\beta^2(y)\beta(w_j)) = (x,\mu_{12}'(\beta(y))),
 \end{eqnarray*}
where the penultimate equality holds since if $\{ w_j^{\ast} \}$ is a dual basis of $\{ w_j \}$, then $\{ \beta(w_j) \}$ is a dual basis of $\{ w_j^{\ast} \}$.

From the self-duality of the Hochschild homology complex (\ref{HH_hom_complex-almostCY}) and $(\mu_{12}')^{\ast} = \mu_{12}' \circ \beta$, we have
\begin{eqnarray*}
HH_i(A)^{\ast} & \cong & HH_{11-i}(A)[3h], \qquad \qquad \; i=1,\ldots,10, \\
HH_{11}(A)^{\ast} & \cong & HH_{12}(A)[6h].
\end{eqnarray*}

The reduced Hochschild homology $\overline{HH}_{\bullet}(A)$ is defined as $\overline{HH}_0(A) = HH_0(A)/S$ and $\overline{HH}_n(A) = HH_n(A)$, $n>0$.

\subsection{The cyclic homology of $A(\mathcal{G},W)$} \label{sect:cyclic_hom}

Before we determine the Hochschild homology of $A(\mathcal{G},W)$ for certain $SU(3)$ $\mathcal{ADE}$ graphs, we introduce cyclic homology.
We begin by introducing the differential graded algebra $\Omega^{\bullet}A$ of non-commutative forms of $A$, and the non-commutative de Rham homology.

The $A$-$A$ bimodule $\Omega^1 A$ of non-commutative relative 1-forms on $A$ is defined as the kernel of the multiplication map $A \otimes_S A \rightarrow A$. The differential graded algebra $\Omega^{\bullet}A$ of non-commutative forms of $A$ is obtained by taking tensor powers of $\Omega^1 A$. The graded commutator in $\Omega^{\bullet}A$ is given by $[\omega,\omega'] = \omega \omega' - (-1)^{|\omega| |\omega'|} \omega' \omega$, where $|\omega| = n$ denotes the homological degree of $\omega \in \Omega^n A$. The reduced non-commutative de Rham homology of $A$ is defined by
$$\overline{H}DR_n (A) := H_n (\Omega^{\bullet}A/(S+[\Omega^{\bullet}A,\Omega^{\bullet}A]),d),$$
where the natural differential $\Omega^{\bullet}A \rightarrow \Omega^{\bullet+1}A$ descends to a de Rham differential on $\Omega^{\bullet}A/(S+[\Omega^{\bullet}A,\Omega^{\bullet}A])$.

Since $A$ is an augmented $S$-algebra, i.e. $A_0 = S$ and there is an augmentation $\varphi:A \rightarrow S$ such that $\varphi(1)=1$, by the non-commutative Poincar\'{e} lemma \cite{kontsevich:1993} (see also \cite[Lemma 4.5]{mejias:2002}), $\overline{H}DR_n (A) = \overline{H}DR_n (S) = 0$ for all $n$. Thus, from \cite[Lemma 3.6.1]{etingof/ginzburg:2007}, there is an exact sequence
\begin{equation} \label{seq:Connes_exact_seq}
0 \longrightarrow \overline{HH}_0(A) \stackrel{B}{\longrightarrow} \overline{HH}_1(A) \stackrel{B}{\longrightarrow} \overline{HH}_2(A) \stackrel{B}{\longrightarrow} \overline{HH}_0(A) \longrightarrow \cdots
\end{equation}
where $B$ is the Connes differential, which is degree-preserving, and the reduced cyclic homology of $A$ can be defined by
$$\overline{HC}_n (A) = \mathrm{ker}(B: \overline{HH}_{n+1}(A) \rightarrow \overline{HH}_{n+2}(A)) = \mathrm{Im}(B: \overline{HH}_{n}(A) \rightarrow \overline{HH}_{n+1}(A)).$$
The usual cyclic homology is related to the reduced cyclic homology by $\overline{HC}_0(A) = HC_0(A)/S$ and $\overline{HC}_n(A) = HC_n(A)$, $n>0$.

The (graded) Euler characteristic of the reduced cyclic homology is the polynomial in $t$ defined by $\chi_{\overline{HC}(A)}(t) = \sum_{i=0}^{\infty} (-1)^i H_{\overline{HC}_i(A)}(t)$.
It turns out to be easier to describe the Euler characteristic of $\mathrm{Sym}(\overline{HC}(A))_+$, where if $\chi_{\overline{HC}(A)}(t) = \sum_{k=0}^{\infty}  a_k t^k$ then $\chi_{\mathrm{Sym}(\overline{HC}(A))_+}(t) = \prod_{k=1}^{\infty} (1-t^k)^{a_k}$.
In \cite[Prop. 3.7.1]{etingof/ginzburg:2007} it was shown that for $A$ the preprojective algebra of a non-Dynkin quiver,
\begin{equation} \label{eqn:Euler=prod_det}
\prod_{k=1}^{\infty} (1-t^k)^{-a_k} = \prod_{s=1}^{\infty} \mathrm{det}H_A(t^s),
\end{equation}
where $H_A(t)$ is the Hilbert series of $A$. The result (\ref{eqn:Euler=prod_det}) was extended to the case where $A$ is a Calabi-Yau algebra of dimension 3 in \cite[Prop. 5.4.9]{ginzburg:2006}.
In the case when $A$ is the almost Calabi-Yau algebra $A=A(\mathcal{G},W)$, the differential graded algebra $\mathfrak{D}_{\bullet} = T_S(V \oplus V^{\ast} \oplus S^{\ast})$ in \cite[Prop. 5.4.9]{ginzburg:2006} is no longer exact. However, we can build a larger free differential graded algebra $\mathfrak{D}_{\bullet}'$ by adding generators $x_n \in \mathfrak{D}_n'$ whose images under the differential give a basis for $H_n(\mathfrak{D}_{\bullet}')$, for each $n>0$. These generators lie in degree $nh$, where $h$ is the Coxeter number of $\mathcal{G}$. Then $\mathfrak{D}_{\bullet}'$ gives a free resolution of $A$, and a correction term corresponding to the numerator $1-Pt^{hs}$ of $H_A(t^s)$ appears in the formula (\ref{eqn:Euler=prod_det}). Thus the result (\ref{eqn:Euler=prod_det}) holds for the almost Calabi-Yau algebra $A=A(\mathcal{G},W)$ (c.f. \cite[Lemma 4.4.1]{etingof/eu:2007} in the case where $A$ is the preprojective algebra of a Dynkin quiver).

\subsection{$HH_0(A)$ for $A=A(\mathcal{G},W)$} \label{sect:HH_0}

In this section we compute the zeroth Hochschild homology $HH_0(A) = \mathrm{ker}(\mu_0')/\mathrm{Im}(\mu_1') = A/[A,A]$ for the simplest graphs, namely the graphs $\mathcal{A}^{(n)}$, $n \geq 4$, $\mathcal{D}^{(3k+3)}$, $k \geq 1$, $\mathcal{A}^{(n)\ast}$, $n \geq 5$, $\mathcal{D}^{(n)\ast}$, $n \geq 5$, $\mathcal{E}^{(8)}$ and $\mathcal{E}^{(8)\ast}$.

For a graded algebra $B = \bigoplus_{k=0}^{\infty} B_k$, let $B_+$ denote the positive degree part $B_+ = \bigoplus_{k=1}^{\infty} B_k$. For any $a,b \in A_+$ such that $r(a) = s(b)$ and $s(a) \neq r(b)$, $[a,b] = ab$, thus any non-cyclic path $ab$ is in $[A,A]$.
For $a,b \in A_+$ such that $r(a) = s(b)$ and $s(a) = r(b)$, $[a,b] = ab-ba$, thus cyclic paths are equivalent in $A/[A,A]$ if one is a cyclic permutation of the other.
Thus to determine $A/[A,A]$ we first consider all cyclic paths in $i A i$ for some $i \in \mathcal{G}_0$, then consider all cyclic paths in $j A j$ which do not pass through the vertex $i \neq j$, for some $j \in \mathcal{G}_0$, and so on.
Note that $S \hookrightarrow A/[A,A]$ since $[i,j]=0$ for all $i,j \in S$ and $[a,b] \subset A_+$ if either $a$ or $b$ have non-zero length.

\subsubsection{The identity $\mathcal{A}^{(n)}$ graphs}

The unique cell system $W$ (up to equivalence) was computed in \cite[Theorem 5.1]{evans/pugh:2009i}.
For the graph $\mathcal{A}^{(n)}$, $n \geq 4$, the space of cyclic paths $(0,0) A_+ (0,0) = 0$. Thus for any vertex $i \neq (0,0)$, any cyclic path $x \in i A_+ i$ which passes through $(0,0)$ is a cyclic permutation of a cyclic path $x' \in (0,0) A_+ (0,0)$.
Similarly, any cyclic path $x \in i A_+ i$ which does not pass through $(0,0)$ can be transformed by a combination the relations in $A$ and cyclic permutations to a cyclic path $x' \in (0,0) A_+ (0,0)$. Thus any cyclic path $x \in i A_+ i$ will be zero in $A/[A,A]$,
and we obtain
\begin{equation}
HH_0(A) \cong S.
\end{equation}

\subsubsection{The orbifold $\mathcal{D}^{(3k+3)}$ graphs}

\begin{figure}[tb]
\begin{center}
  \includegraphics[width=55mm]{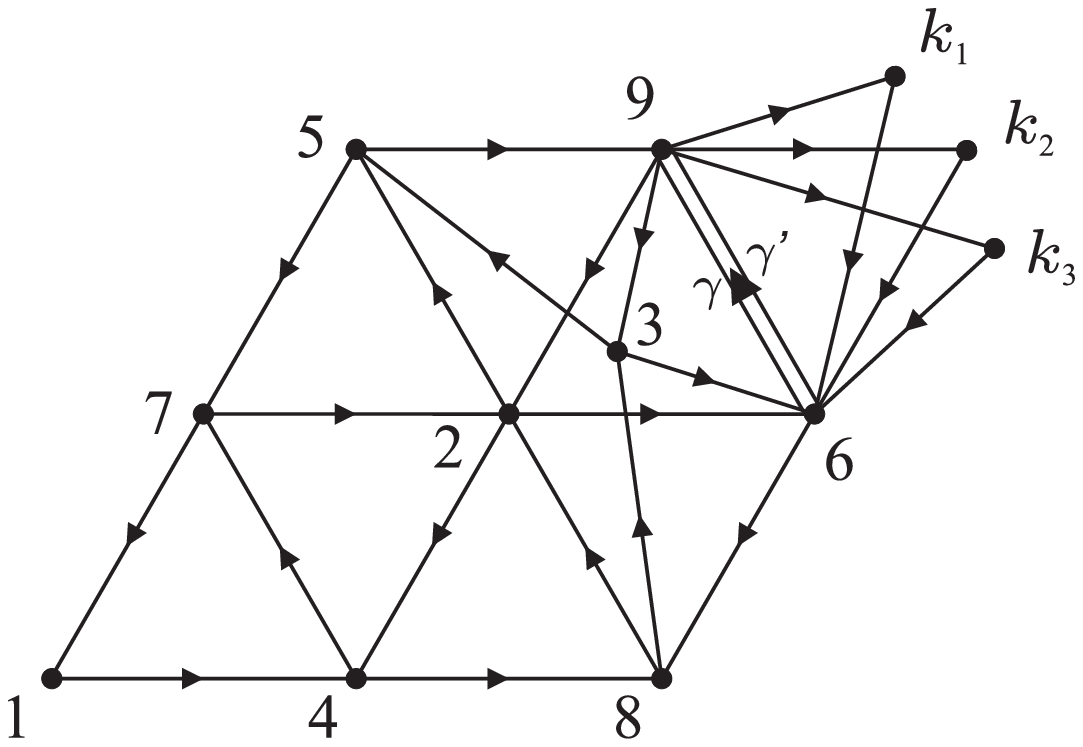}\\
  \caption{Graph $\mathcal{D}^{(9)}$} \label{fig:D(9)}
\end{center}
\end{figure}

We now consider the graphs $\mathcal{D}^{(3k+3)}$, $k \geq 1$, which are $\mathbb{Z}_3$-orbifolds of $\mathcal{A}^{(3k+3)}$. The graph $\mathcal{D}^{(9)}$ is illustrated in Figure \ref{fig:D(9)}. The weights $W(\triangle)$ for $\mathcal{A}^{(3k+3)}$ are invariant under the $\mathbb{Z}_3$ symmetry of the graph given by rotation by $2\pi/3$. Thus there is an orbifold solution for the cell system $W$ on $\mathcal{D}^{(3k+3)}$ where the weights $W(\triangle)$ are given by the corresponding weights for $\mathcal{A}^{(3k+3)}$ \cite[Theorem 6.2]{evans/pugh:2009i}. More precisely, excluding triangles $\triangle$ which contain one of the triplicated vertices $(k,k)_l$, the weight $W(\triangle_{i_1,i_2,i_3})$ for the triangle $\triangle_{i_1,i_2,i_3} = i_1 \rightarrow i_2 \rightarrow i_3 \rightarrow i_1$ on $\mathcal{D}^{(3k+3)}$ is given by the weight $W(\triangle_{i_1^{(0)},i_2^{(1)},i_3^{(2)}}) = W(\triangle_{i_1^{(1)},i_2^{(2)},i_3^{(0)}}) = W(\triangle_{i_1^{(2)},i_2^{(0)},i_3^{(1)}})$ for $\mathcal{A}^{(3k+3)}$, where $i_l^{(0)}$, $i_l^{(1)}$, $i_l^{(2)}$ are the three vertices of $\mathcal{A}^{(3k+3)}$ which are identified under the $\mathbb{Z}_3$ action to give the vertex $i_l$ of $\mathcal{D}^{(3k+3)}$, $l=1,2,3$.
If for a triangle $\triangle_{i_1,i_2,i_3}$ on $\mathcal{D}^{(3k+3)}$ there is no choice of vertices $i_1^{(j_1)}$, $i_2^{(j_2)}$, $i_3^{(j_3)}$ on $\mathcal{A}^{(3k+3)}$ which lie on a closed loop of length three $i_1^{(j_1)} \rightarrow i_2^{(j_2)} \rightarrow i_3^{(j_3)} \rightarrow i_1^{(j_1)}$, then we have $W(\triangle_{i_1,i_2,i_3}) = 0$.
The weight $W(\triangle)$ for a triangle $\triangle$ which contain one of the triplicated vertices $(k,k)_l$ is just given by one third of the weight for the corresponding triangle on $\mathcal{A}^{(3k+3)}$.
Thus the relations for $\mathcal{D}^{(3k+3)}$ are given precisely by the relations for $\mathcal{A}^{(3k+3)}$, except for the relations $\rho_{\gamma}$, $\rho_{\gamma'}$, which involve the triplicated vertices $(k,k)_l$.

Any cyclic path on $\mathcal{A}^{(3k+3)}$ yields a cyclic paths on $\mathcal{D}^{(3k+3)}$ by the above orbifold procedure. These cyclic paths will be zero in $A/[A,A]$, except for those which pass along the double edge of $\mathcal{D}^{(3k+3)}$ -- although these paths can be made to pass through $(0,0)$ in $A'/[A',A']$ for $A'=A(\mathcal{A}^{(3k+3)},W)$, when we do this for $A=A(\mathcal{D}^{(3k+3)},W)$ we obtain a cyclic path which passes through the vertex 1 of $\mathcal{D}^{(3k+3)}$ which corresponds to $(0,0)$ on $\mathcal{A}^{(3k+3)}$, but also a linear combination of cyclic paths which do not pass through the vertex 1, due to the fact that relations involving the double edge are not of the form $x=\lambda x'$ for basis paths $x,x' \in A$.
There are also cyclic paths in $A$ which do not come from cyclic paths in $A'$ by the orbifold procedure. These paths must necessarily pass along the double edge $(\gamma,\gamma')$ of $\mathcal{D}^{(3k+3)}$. Using the relations in $A$ and cyclic permutations, we can transform any such cyclic path, necessarily of length $3j$, $j \in \mathbb{N}$, due to the three-colourability of $\mathcal{D}^{(3k+3)}$, to a linear combination of cyclic basis paths $[(i_1 i_2 k_l i_1)^j]$, $l=1,2,3$, where $i_1 = s(\gamma)$, $i_2 = r(\gamma)$, $k_l := (k,k)_l$, and $x^m$ denotes the path $x x \cdots x$ ($m$ times). These basis paths are not equivalent in $[A,A]$, except when $j=k$ where $[(i_1 i_2 k_1 i_1)^{k}] = [(i_1 i_2 k_2 i_1)^{k}] = [(i_1 i_2 k_3 i_1)^{k}]$, thus
\begin{equation}
HH_0(A) \cong S \oplus C,
\end{equation}
where the graded vector space $C = \bigoplus_{j=1}^{k-1} \bigoplus_{l=1}^3 \mathbb{C} [(i_1 i_2 k_l i_1)^j] \oplus \mathbb{C} [(i_1 i_2 k_1 i_1)^{k}]$, and has Hilbert series $H_C(t) = \sum_{j=1}^{k-1} 3t^{3j} + t^{3k}$.

\subsubsection{The conjugate $\mathcal{A}^{\ast}$ graphs}

\begin{figure}[tb]
\begin{minipage}[t]{6cm}
\begin{center}
  \includegraphics[width=35mm]{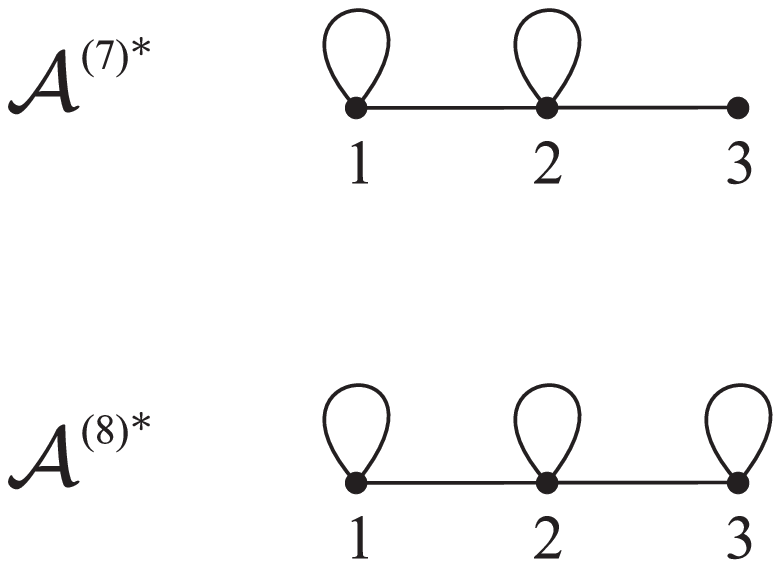}\\
  \caption{Graphs $\mathcal{A}^{(7)\ast}$, $\mathcal{A}^{(8)\ast}$} \label{fig:A(7,8)star}
\end{center}
\end{minipage}
\hfill
\begin{minipage}[t]{9cm}
\begin{center}
  \includegraphics[width=90mm]{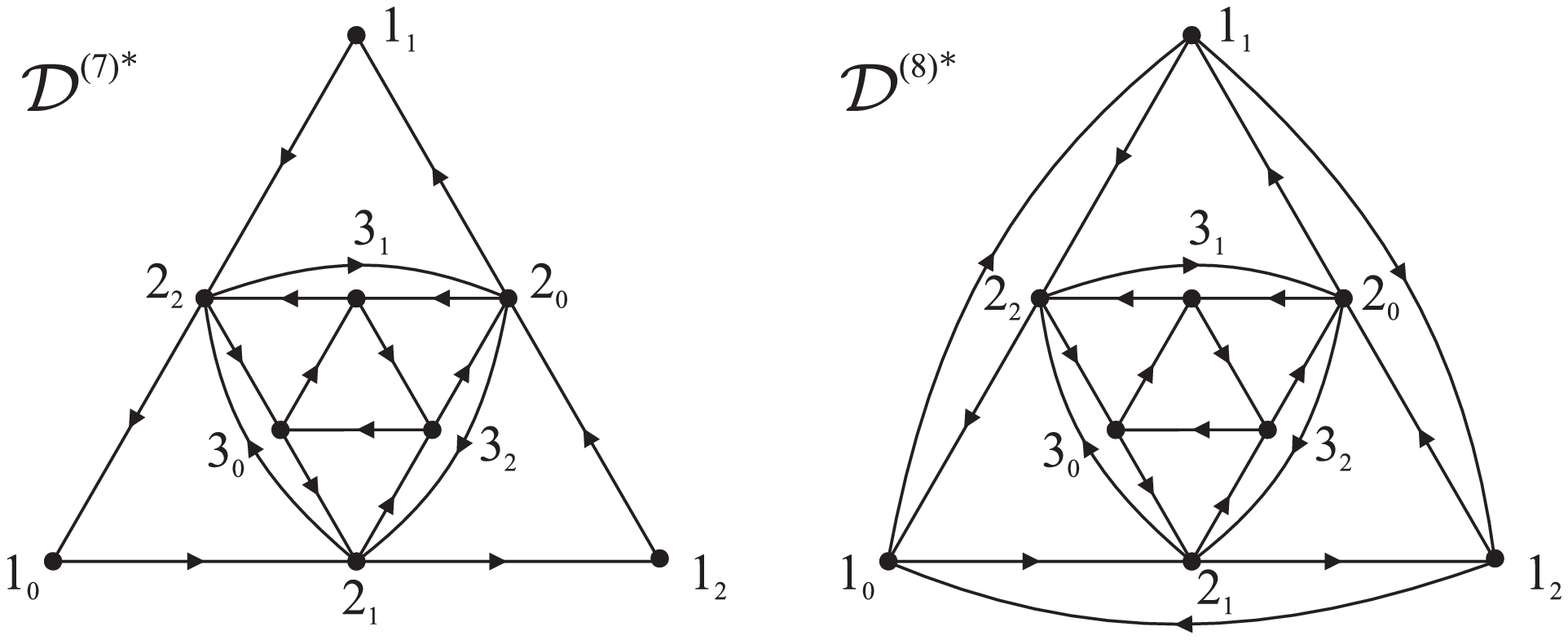}\\
  \caption{Graphs $\mathcal{D}^{(7)\ast}$, $\mathcal{D}^{(8)\ast}$} \label{fig:D(7,8)star}
\end{center}
\end{minipage}
\end{figure}

The unique cell system $W$ (up to equivalence) was computed in \cite[Theorems 7.1, 7.3 \& 7.4]{evans/pugh:2009i}, and we use the same notation for the cells here.
The $\mathcal{A}^{(n)\ast}$ graphs are illustrated in \cite[Figure 11]{evans/pugh:2009i}. We illustrate the cases $n=7,8$ here in Figure \ref{fig:A(7,8)star}. The numbering of the vertices of $\mathcal{A}^{(2m+1)\ast}$ that we use here is the same as that in \cite{evans/pugh:2010ii}, but the reverse of that used in \cite{evans/pugh:2009i}.
The relations in $A(\mathcal{A}^{(n)\ast},W)$ are
\begin{eqnarray}
& W_{112}[121] + W_{111}[111] = 0, & \nonumber \label{eqn:rel-Astar-i} \\
& W_{p-1,p,p}[p(p-1)j] + W_{p,p,p}[ppp] + W_{p,p,p+1}[p(p+1)p] = 0, & \label{eqn:rel-Astar-ii}  \\
& W_{p,p,p+1}[pp(p+1)] + W_{p,p+1,p+1}[p(p+1)(p+1)] = 0, & \label{eqn:rel-Astar-iii} \\
& W_{p,p,p+1}[(p+1)pp] + W_{p,p+1,p+1}[(p+1)(p+1)p] = 0, & \label{eqn:rel-Astar-iv}
\end{eqnarray}
where $p=2,\ldots,r-1$ in (\ref{eqn:rel-Astar-ii}), and $p=1,\ldots,p'$ in (\ref{eqn:rel-Astar-iii}), (\ref{eqn:rel-Astar-iv}), where $r = \lfloor (n-1)/2 \rfloor$, $p'=r-1$ for even $n$, and $p'=r-2$ for odd $n$.
For even $n$ we have the extra relation $W_{r-1,r,r}[r(r-1)r] + W_{r,r,r}[rrr] = 0$, and for odd $n$ we have the extra relation $[r(r-1)(r-1)] = [(r-1)(r-1)r] = 0$.

We first consider the even case $n=2m+2$.
Clearly all loops $[pp]$ of length 1 are in $A/[A,A]$, $p=1,\ldots,m$.
Let $d_l^p := \mathrm{dim}(p A_l p)$. From the Hilbert series for $A$, we see that $d_{2k}^p = d_{2k+1}^p = p = d_{2m-2k-1}^p = d_{2m-2k-2}^p$ if $p \leq k$ or $p \geq m-k+1$, and $d_{2k}^p = d_{2k+1}^p = k+1 = d_{2m-2k-1}^p = d_{2m-2k-2}^p$ if $k+1 \leq p \leq m-k$, for $2k \leq m-1$. Then $\mathrm{dim}(A^S_{2k}) = \mathrm{dim}(A^S_{2k+1}) = (k+1)(m-k) = \mathrm{dim}(A^S_{2m-2k-1}) = \mathrm{dim}(A^S_{2m-2k-2})$ for $2k \leq m-1$.

Each commutator of the form $[[l(l+1)],[(l+1)l]] = [l(l+1)l] - [(l+1)l(l+1)]$ yields a relation between linearly independent paths of length 2 in $A/[A,A]$, $l=1,\ldots,m-1$. There are $m-1$ such relations, thus the dimension of $(A/[A,A])_2$ is $2(m-1) - (m-1) = m-1$. Similarly the dimension of $(A/[A,A])_3$ is $m-1$.
Each commutator of the form $C_p = [[(p-1)ppp],[p(p-1)]]$, $p=2,\ldots,m$, and $C_p' = [[(p-1)p(p+1)],[(p+1)p(p-1)]]$, $p=2,\ldots,m-1$, yield relations between linearly independent paths of length 4 in $A/[A,A]$. There is one basis path in $1 A^S_4 1$, which we may take to be $[11111]$. Let $w_1$ denote the basis element given by its image in $A/[A,A]$. Since $d_4^2 = 2$, the dimension of $2 A^S_4 2$ is 2. However, the basis can be chosen such that one of the basis paths is identified with $w_1$ in $A/[A,A]$ by $C_2$, thus we obtain one new basis path $w_2 \in (A/[A,A])_4$, which may be chosen to be $[22222]$. Similarly, the dimension of $p A^S_4 p$ is 3, $p=3,\ldots,m-2$, and the basis can be chosen such that two of the basis paths are identified with linear combinations of $w_1, w_2, \ldots, w_{p-1}$ in $A/[A,A]$ by $C_p$ and $C_p'$. Thus we obtain one new basis path $w_p \in (A/[A,A])_4$ for each $p=3,\ldots,m-2$, which may be chosen to be $[ppppp]$. The dimension of $(m-1) A^S_4 (m-1)$ is 2, but by $C_{m-1}$, $C_{m-1}'$ any such path can be identified with a linear combination of $w_1, w_2, \ldots, w_{m-2}$ in $A/[A,A]$. Similarly the single basis path in $m A^S_4 m$ can be identified with a linear combination of $w_1, w_2, \ldots, w_{m-2}$ in $A/[A,A]$. Thus we obtain a basis $\{ w_1, \ldots, w_{m-2} \}$ for $(A/[A,A])_4$.
By a similar argument, we see that the dimension of $(A/[A,A])_k$ is $m- \lfloor k/2 \rfloor$ for all $k=0,1,\ldots,2m-1$, with basis paths $[ppp \cdots p]$, for $p=1,\ldots,m- \lfloor k/2 \rfloor$.
Thus
\begin{equation}
HH_0(A) \cong S \oplus C,
\end{equation}
where the graded vector space $C$ has Hilbert series $H_C(t) = \sum_{j=1}^{2m-1} (m- \lfloor j/2 \rfloor) t^{j}$.

We now consider the odd case $n=2m+1$.
Again, all loops $[pp]$ of length 1 are in $A/[A,A]$, this time for $p=1,\ldots,m-1$ (note that there is no edge from vertex $m$ to $m$ on $\mathcal{A}^{(2m+1)\ast}$).
From the Hilbert series for $A$, we see that $d_{2k}$ is given by the same formula as for the even case $n=2m+2$, for $2k \leq m-1$, whilst $d_{2k+1}^p = p = d_{2m-2k-2}^p$ if $p \leq k$ or $p \geq m-k$, $d_{2k+1}^p = k+1 = d_{2m-2k-3}^p$ if $k+1 \leq i \leq m-k-1$, and $d_{2k+1}^m = 0$, for $2k \leq m-2$.
Then $\mathrm{dim}(A^S_{2k}) = (k+1)(m-k) = \mathrm{dim}(A^S_{2m-2k-2})$ for $2k \leq m-1$, and $\mathrm{dim}(A^S_{2k-1}) = (k+1)(m-k-1) = \mathrm{dim}(A^S_{2m-2k-3})$ for $2k \leq m-2$.
By a similar argument as for the even case above, we see that the dimension of $(A/[A,A])_k$ is $m- \lfloor (k+1)/2 \rfloor$ for all $k=0,1,\ldots,2m-2$, with basis paths $[ppp \cdots p]$, for $p=1,\ldots,m- \lfloor (k+1)/2 \rfloor$.
Thus
\begin{equation}
HH_0(A) \cong S \oplus C,
\end{equation}
where the graded vector space $C$ has Hilbert series $H_C(t) = \sum_{j=1}^{2m-2} (m- \lfloor (j+1)/2 \rfloor) t^{j}$.

\subsubsection{The conjugate orbifold $\mathcal{D}^{\ast}$ graphs}

The graphs $\mathcal{D}^{(n)\ast}$ are (three-colourable) unfolded versions, or $\mathbb{Z}_3$-orbifolds, of the graphs $\mathcal{A}^{(n)\ast}$, where we replace every vertex $v$ of $\mathcal{A}^{(n)\ast}$ by three vertices $v_0$, $v_1$, $v_2$, where $v_a$ is of colour $a$, such that there are edges $v_0 \rightarrow w_1$, $v_1 \rightarrow w_2$ and $v_2 \rightarrow w_0$ if and only if there is an edge $v \rightarrow w$ on $\mathcal{A}^{(n)\ast}$. The graphs $\mathcal{D}^{(7)\ast}$, $\mathcal{D}^{(8)\ast}$ are illustrated in Figure \ref{fig:D(7,8)star}.

Due to the three-colourability of the graph $\mathcal{D}^{(n)\ast}$, a closed loop on $\mathcal{A}^{(n)\ast}$ will only be a closed loop on $\mathcal{D}^{(n)\ast}$ if it has length $3k$, $k \geq 0$, and for each such closed loop on $\mathcal{A}^{(n)\ast}$, there are three corresponding closed loops of length $3k$ on $\mathcal{D}^{(n)\ast}$. However, these three closed loops are identified in $A/[A,A]$, which can be seen as follows. As in the case of $\mathcal{A}^{(n)\ast}$, $(A/[A,A])_{3k}$ is generated by paths of the form $[p_l p_{l+1} p_{l+2} p_l \cdots p_{l+3k}]$, for $l=0,1,2 \textrm{ mod } 3$ and $p=1,\ldots,r$, where $r = m- \lfloor 3k/2 \rfloor$ for $n=2m+2$ and $r = m- \lfloor (3k+1)/2 \rfloor$ for $n=2m+1$. Since $[p_l p_{l+1} p_{l+2} p_l \cdots p_{l}]$ is a cyclic permutation of $[p_{l+1} p_{l+2} p_l p_{l+1} \cdots p_{l+1}]$, we see that for $l=0,1,2 \textrm{ mod } 3$, the cyclic paths $[p_l p_{l+1} p_{l+2} p_l \cdots p_{l+3k}]$ are identified in $A/[A,A]$. Thus $(A/[A,A])_{3k}$ has a basis given by $[p_0 p_1 p_2 p_0 \cdots p_0]$, for $p=1,\ldots,r$, and $(A/[A,A])_{k'} = 0$ for $k' \not \equiv 0 \textrm{ mod } 3$.
Then
\begin{equation}
HH_0(A) \cong S \oplus C,
\end{equation}
where for $n=2m+2$ the graded vector space $C$ has Hilbert series $H_C(t) = \sum_{j=1}^{\lfloor (2m-1)/3 \rfloor} (m- \lfloor 3j/2 \rfloor) t^{3j}$,
whilst for $n=2m+1$, $H_C(t) = \sum_{j=1}^{\lfloor (2m-2)/3 \rfloor} (m- \lfloor (3j+1)/2 \rfloor) t^{3j}$.

\subsubsection{The graph $\mathcal{E}^{(8)}$ for the conformal embedding $SU(3)_5 \subset SU(6)_1$}

\begin{figure}[tb]
\begin{minipage}[t]{7.5cm}
\begin{center}
  \includegraphics[width=35mm]{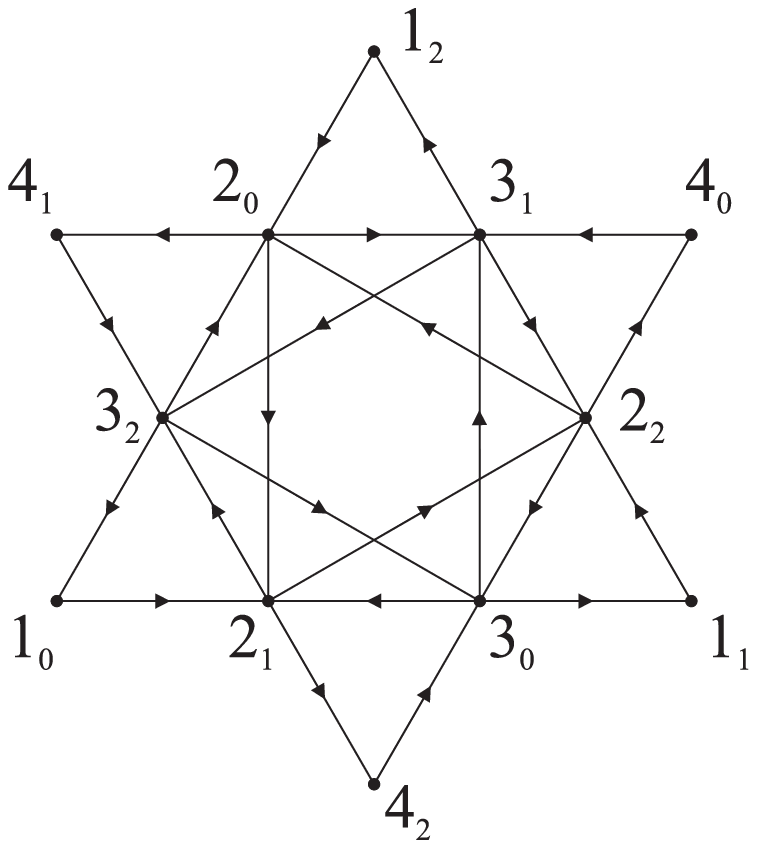}\\
  \caption{Graph $\mathcal{E}^{(8)}$} \label{fig:E(8)}
\end{center}
\end{minipage}
\hfill
\begin{minipage}[t]{7.5cm}
\begin{center}
  \includegraphics[width=15mm]{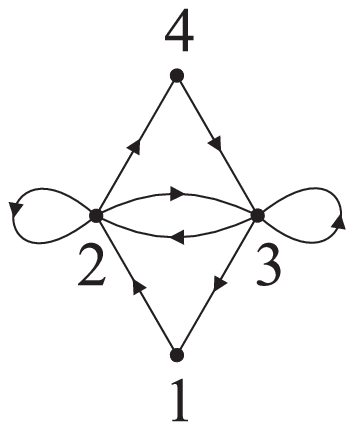}\\
  \caption{Graph $\mathcal{E}^{(8)\ast}$} \label{fig:E(8)star}
\end{center}
\end{minipage}
\end{figure}

We now consider the graph $\mathcal{E}^{(8)}$, illustrated in Figure \ref{fig:E(8)}. The unique cell system $W$ (up to equivalence) was computed in \cite[Theorem 9.1]{evans/pugh:2009i}. The quotient algebra $A$ has the relations, for $l=0,1,2$,
\begin{eqnarray*}
& [1_l 2_{l+1} 3_{l+2}] = [2_{l+1} 3_{l+2} 1_l] = [3_l 2_{l+1} 4_{l+2}] = [4_{l+2} 3_l 2_{l+1}] = 0, & \\
& \sqrt{[3]} [2_l 2_{l+1} 2_{l+2}] = -[2_l 3_{l+1} 2_{l+2}], \qquad \sqrt{[3]}[3_l 3_{l+1} 3_{l+2}] = [3_l 2_{l+1} 3_{l+2}], & \\
& \frac{-\sqrt{[3]}}{[2]} [3_l 1_{l+1} 2_{l+2}] = [3_l 2_{l+1} 2_{l+2}] + [3_l 3_{l+1} 2_{l+2}], & \\
& \frac{-\sqrt{[3]}}{[2]} [2_l 4_{l+1} 3_{l+2}] = [2_l 2_{l+1} 3_{l+2}] + [2_l 3_{l+1} 3_{l+2}]. &
\end{eqnarray*}

The only cyclic paths in $A_+$ are of the form $[2_l 2_{l+1} 2_{l+2} 2_l]$, $[3_l 3_{l+1} 3_{l+2} 3_l]$, $l=1,\ldots,3$. Now $[2_l 2_{l+1} 2_{l+2} 2_l] = [2_{l+1} 2_{l+2} 2_l 2_{l+1}]$ by cyclic permutation, and $\sqrt{[3]} [2_l 2_{l+1} 2_{l+2} 2_l] = -[2_l 3_{l+1} 2_{l+2} 2_l] = [3_{l+1} 2_{l+2} 2_l 3_{l+1}] = -[3_{l+1} 2_{l+2} 3_l 3_{l+1}] - (\sqrt{[2]}/\sqrt{[4]}) [3_{l+1} 2_{l+2} 4_l 3_{l+1}] = -\sqrt{[3]} [3_{l+1} 3_{l+2} 3_l 3_{l+1}]$ in $A/[A,A]$, where the second and last equalities follows by cyclic permutation and the others follow from the relations in $A$. Thus we see that all cyclic paths in $A_+$ are identified in $A/[A,A]$ so that
\begin{equation}
HH_0(A) \cong S \oplus \mathbb{C} [2_0 2_1 2_2 2_0].
\end{equation}

\subsubsection{The graph $\mathcal{E}^{(8)\ast}$ for the orbifold of the conformal embedding $SU(3)_5 \subset SU(6)_1 \rtimes \mathbb{Z}_3$}

Consider the graph $\mathcal{E}^{(8)\ast}$, illustrated in Figure \ref{fig:E(8)star}. The unique cell system $W$ (up to equivalence) was computed in \cite[Theorem 10.1]{evans/pugh:2009i}. The quotient algebra $A$ has the relations
\begin{eqnarray*}
& [123] = [231] = [324] = [432] = 0, \qquad \qquad [222] = \frac{-1}{\sqrt{[3]}} [232], \qquad \qquad [333] = \frac{1}{\sqrt{[3]}} [323], & \\
& \frac{-\sqrt{[3]}}{[2]} [312] = [322] + [332], \qquad \qquad \frac{-\sqrt{[3]}}{[2]} [243] = [223] + [233]. &
\end{eqnarray*}

Clearly the single edges $[22]$, $[33]$ are not in $[A,A]$, since the relations in $A$ only change paths of length $>1$, and edges are invariant under cyclic permutation. We have the relation $\sqrt{[3]} [222 \cdots 2] = -[232 \cdots 2] = -[322 \cdots 23]$ for paths of length $r$ in $A/[A,A]$, $2 \leq r \leq 5$, where the first equality follows from the relation in $A$ and the second follows by cyclic permutation.
Thus in $A/[A,A]$, for $r=2$, we obtain $\sqrt{[3]} [222] = -[323] = -\sqrt{[3]} [333]$ by the relations in $A$.
For $r=3$, we have $\sqrt{[3]} [2222] = -[3223] = [3323] +([2]/\sqrt{[3]})[3123] = \sqrt{[3]}[3333]$, by the relations in $A$, since the subpath $[123] = 0$ in $A$.
For $r=4$, we have $[3] [22222] = -\sqrt{[3]}[32223] = [32323] = \sqrt{[3]}[33323] = [3] [33333]$, by the relations in $A$, but also $\sqrt{[3]} [22222] = -[32223] = [33223] +([2]/\sqrt{[3]})[31223] = -[33233] -([2]/\sqrt{[3]})[33243] +([2]/\sqrt{[3]})[23122] = -\sqrt{[3]}[33333]$, by the relations in $A$, since the subpaths $[324] = 0 = [123]$ in $A$, and we have used the cyclic permutation relation in the penultimate equality. Then in $A/[A,A]$, we see that $[22222] = 0 = [33333]$.
For $r=5$ we have $[3] [222222] = -\sqrt{[3]}[322223] = [323223] = \sqrt{[3]}[333223] = -[3][333333]$ in $A/[A,A]$, by the relations in $A$.
Thus
\begin{equation}
HH_0(A) \cong S \oplus C,
\end{equation}
where the graded vector space $C = \mathbb{C} \{ [22], [33] \} \oplus  \mathbb{C} [222] \oplus \mathbb{C} [2222] \oplus \mathbb{C} [222222]$, and has Hilbert series $H_C(t) = 2t + t^2 + t^3 + t^5$.

\subsection{Determining the Hochschild homology of $A(\mathcal{G},W)$ for trivial Nakayama automorphism} \label{Sect:HH(A)-trivial_beta}

In this section we determine the Hochschild and cyclic homology for the graphs $\mathcal{D}^{(3k)}$, $k \geq 2$, $\mathcal{A}^{(n)\ast}$, $n \geq 4$, $\mathcal{D}^{(3k)\ast}$, $k \geq 2$, and $\mathcal{E}^{(8)\ast}$.
Here the almost Calabi-Yau algebra $A$ has trivial Nakayama automorphism.

In this case, the Hochschild homology of $A$ has minimal period at most 4, thus we have $HH_i(A)^{\ast}[h] \cong HH_{3-i}(A)$, $i=1,2$, and $HH_i(A)^{\ast} \cong HH_{7-i}(A)$, $i=3,4$.
From the exactness of (\ref{seq:Connes_exact_seq}) we see that $\mathrm{ker}(B: \overline{HH}_0(A) \rightarrow \overline{HH}_1(A)) = 0$, and since the Connes differential $B$ preserves degrees, we have $\overline{HH}_1(A) \cong C \oplus X$, for some graded vector space $X$ which lives in degrees 1 to $h-2$. Then $\overline{HH}_2(A) \cong C^{\ast}[h] \oplus X^{\ast}[h]$, where $C \cong \overline{HH}_0(A)$ and $X^{\ast}[h]$ lives in degrees 2 to $h-1$. Now $B:\overline{HH}_1(A) \rightarrow \overline{HH}_2(A)$ restricts to an isomorphism $X \stackrel{\cong}{\longrightarrow} X^{\ast}[h]$ since (\ref{seq:Connes_exact_seq}) is exact, and since it preserves degrees, $X$ only lives in degrees 2 to $h-2$. A similar argument shows that $\mathrm{ker}(B: \overline{HH}_0(A) \rightarrow \overline{HH}_1(A)) \cong X^{\ast}[h]$, so that $\overline{HH}_3(A) \cong C^{\ast}[h] \oplus K'$, where the graded vector space $K'$ lives in degrees 3 to $h$, and $\overline{HH}_4(A) \cong C[h] \oplus K'^{\ast}[h]$, where $K'^{\ast}[h]$ lives in degrees $h$ to $2h-3$. Since $B:\overline{HH}_3(A) \rightarrow \overline{HH}_4(A)$ restricts to an isomorphism $K' \stackrel{\cong}{\longrightarrow} K'^{\ast}[h]$, we see that $K'$ lives only in degree $h$. We will write $K'=K[h]$ where $K$ is a vector space which lives in degree 0, so that $\overline{HH}_3(A) \cong C^{\ast}[h] \oplus K[h]$ and $\overline{HH}_4(A) \cong C[h] \oplus K^{\ast}[h]$.

Thus for any almost Calabi-Yau algebra $A$ with trivial Nakayama automorphism
$$\begin{array}{cccccccc}
&& 0 &&&& \qquad & \\
(\textrm{min deg},\textrm{max deg}) && \downarrow &&&&& \\
(0,h-3) && \overline{HH}_0(A) & \cong & C &&& \overline{HC}_0(A) \cong C \\
&& \hspace{-3mm} {\scriptstyle B} \downarrow && \hspace{-2mm} {\scriptstyle \cong} \downarrow &&& \\
(1,h-2) && \overline{HH}_1(A) & \cong & C & \hspace{-3mm} \oplus X && \overline{HC}_1(A) \cong X \\
&& \hspace{-3mm} {\scriptstyle B} \downarrow &&& \hspace{-3mm} {\scriptstyle \cong} \downarrow && \\
(2,h-1) && \overline{HH}_2(A) & \cong & C^{\ast}[h] & \hspace{-3mm} \oplus X^{\ast}[h] && \overline{HC}_2(A) \cong C^{\ast}[h] \\
&& \hspace{-3mm} {\scriptstyle B} \downarrow && \hspace{-2mm} {\scriptstyle \cong} \downarrow &&& \\
(3,h) && \overline{HH}_3(A) & \cong & C^{\ast}[h] & \hspace{-3mm} \oplus K[h] && \overline{HC}_3(A) \cong K[h] \\
&& \hspace{-3mm} {\scriptstyle B} \downarrow &&& \hspace{-3mm} {\scriptstyle \cong} \downarrow && \\
(h,2h-3) && \overline{HH}_4(A) & \cong & C[h] & \hspace{-3mm} \oplus K^{\ast}[h] && \overline{HC}_4(A) \cong C[h] \\
&& \downarrow &&&&& \\
&& \vdots &&&&& \\
\end{array}$$
where $X$ lives in degrees 2 to $h-2$, $K$ lives in degree 0, and $\overline{HH}_{5+i}(A) \cong \overline{HH}_{1+i}(A)[h]$, $\overline{HC}_{4+i}(A) \cong \overline{HC}_i(A)[h]$ for $i \geq 0$. Since $C$ is known, $X$ and $K$ can be determined from the Euler characteristic $\chi_{\overline{HC}(A)}(t)$ as they live in different degrees.

\subsubsection{The graphs $\mathcal{D}^{(3k)}$}

We consider the cases $\mathcal{D}^{(6k)}$, $\mathcal{D}^{(6k+3)}$ separately, $k \geq 1$.
For the graph $\mathcal{D}^{(6k)}$, $k \geq 1$, $\mathrm{det}(H_A(t)) = (1-t^{6k})^{4(3k(k-1)+2)/2}(1-t^{3k})/(1-t^3)^3$, thus $\chi_{\overline{HC}(A)}(t) = (\sum_{j=1}^{2k-1} 3t^{3j} - t^{3k} - (6k(k-1)+2)t^{6k})/(1-t^{6k})$.
Then since $C$ has Hilbert series $H_C(t) = \sum_{j=1}^{2k-2} 3t^{3j} + t^{6k-3}$, we see that $H_X(t) = t^3 + \sum_{j=2}^{2k-2} 3t^{3j} + t^{3k} + t^{6k-3}$ and $H_K(t) = 6k(k-1)+2$, and we obtain:

\begin{Thm}
The Hochschild and cyclic homology of $A=A(\mathcal{D}^{(6k)},W)$, $k \geq 1$, where $W$ is equivalent to one of the cell systems constructed in \cite{evans/pugh:2009i}, is given by
$$\begin{array}{lcl}
HH_0(A) \cong S \oplus C, & \qquad & HC_0(A) \cong S \oplus C, \\
HH_1(A) \cong C \oplus X, && HC_1(A) \cong X, \\
HH_2(A) \cong C^{\ast}[6k] \oplus X^{\ast}[6k], && HC_2(A) \cong C^{\ast}[6k], \\
HH_3(A) \cong C^{\ast}[6k] \oplus K[6k], && HC_3(A) \cong K[6k], \\
HH_4(A) \cong C[6k] \oplus K^{\ast}[6k], && HC_4(A) \cong C[6k], \\
HH_{4+i}(A) \cong HH_i(A)[6k], \quad i \geq 1, && HC_{4+i}(A) \cong HC_i(A)[6k], \quad i \geq 1,
\end{array}$$
where the graded vector spaces $C$, $X$, $K$ have Hilbert series $H_C(t) = \sum_{j=1}^{2k-2} 3t^{3j} + t^{6k-3}$, $H_X(t) = t^3 + \sum_{j=2}^{2k-2} 3t^{3j} + t^{3k} + t^{6k-3}$ and $H_K(t)=6k(k-1)+2$ respectively, where for $k=1$, $H_X(t) = 0$.
\end{Thm}

For $\mathcal{D}^{(6k+3)}$, $k \geq 1$, $\mathrm{det}(H_A(t)) = (1-t^{6k+3})^{6k^2+3}/(1-t^3)^3$, thus $\chi_{\overline{HC}(A)}(t) = (\sum_{j=1}^{2k} 3t^{3j} - 6k^2 t^{6k+3})/(1-t^{6k+3})$. Since $C$ has Hilbert series $H_C(t) = \sum_{j=1}^{2k-1} 3t^{3j} + t^{6k}$, we see that $H_X(t) = t^3 + \sum_{j=2}^{2k-1} 3t^{3j} + t^{6k}$ and $H_K(t) = 6k^2$, and we obtain:

\begin{Thm}
The Hochschild and cyclic homology of $A=A(\mathcal{D}^{(6k+3)},W)$, $k \geq 1$, where $W$ is equivalent to one of the cell systems constructed in \cite{evans/pugh:2009i}, is given by
$$\begin{array}{lcl}
HH_0(A) \cong S \oplus C, & \qquad & HC_0(A) \cong S \oplus C, \\
HH_1(A) \cong C \oplus X, && HC_1(A) \cong X, \\
HH_2(A) \cong C^{\ast}[6k+3] \oplus X^{\ast}[6k+3], && HC_2(A) \cong C^{\ast}[6k+3], \\
HH_3(A) \cong C^{\ast}[6k+3] \oplus K[6k+3], && HC_3(A) \cong K[6k+3], \\
HH_4(A) \cong C[6k+3] \oplus K^{\ast}[6k+3], && HC_4(A) \cong C[6k+3], \\
HH_{4+i}(A) \cong HH_i(A)[6k+3], \quad i \geq 1, && HC_{4+i}(A) \cong HC_i(A)[6k+3], \quad i \geq 1,
\end{array}$$
where the graded vector spaces $C$, $X$, $K$ have Hilbert series $H_C(t) = \sum_{j=1}^{2k-1} 3t^{3j} + t^{6k}$, $H_X(t) = t^3 + \sum_{j=2}^{2k-1} 3t^{3j} + t^{6k}$ and $H_K(t)=6k^2$ respectively.
\end{Thm}

\subsubsection{The $\mathcal{A}^{\ast}$ graphs}

Let $D_m$, $D_m'$ denote the determinant of the denominator $1 - \Delta_{\mathcal{G}} t + \Delta_{\mathcal{G}}^T t^2 - t^3$ of $H_A(t)$ for $\mathcal{G} = \mathcal{A}^{(2m+2)\ast}, \mathcal{A}^{(2m+1)\ast}$ respectively, where $\Delta_{\mathcal{G}}$ denotes the adjacency matrix of $\mathcal{G}$, and let $T_1 = 1-t+t^2-t^3$, $T_2 = t^2-t$.
From the properties of determinants we can deduce the recursion relations $D_m = T_1 D_{m-1} - T_2^2 D_{m-2}$ and $D_m' = (1-t^3) D_{m-1} - T_2^2 D_{m-2}$, for $m \geq 3$, and $D_1 = T_1$, $D_2 = T_1^2 - T_2^2$. It is easy to show by induction on $m$ that $D_m = (1-t)^m(1-t^{2m+2})/(1-t^2)$, and thus $D_m' = (1-t)^{m-1}(1-t^{2m+1})$.
Then for $\mathcal{A}^{(2m+2)\ast}$, $\mathrm{det}H_A(t) = (1-t^{2m+2})^m D_m^{-1} = (1-t^2)(1-t^{2m+2})^{m-1}/(1-t)^m$, thus $\chi_{\overline{HC}(A)}(t) = (mt+(m-1)t^2+mt^3+(m-1)t^4+\cdots+(m-1)t^{2m}+mt^{2m+1})/(1-t^{2m+2})$. Then $H_X(t) = 0 = H_K(t)$, i.e. $X=0=K$.
For $\mathcal{A}^{(2m+1)\ast}$, $\mathrm{det}H_A(t) = (1-t^{2m+1})^m (D_m')^{-1} = (1-t^{2m+1})^{m-1}/(1-t)^{m-1}$, thus $\chi_{\overline{HC}(A)}(t) = ((m-1)t+(m-1)t^2+(m-1)t^3+\cdots+(m-1)t^{2m})/(1-t^{2m+1})$. Then we again deduce that $H_X(t) = 0 = H_K(t)$, and we obtain:

\begin{Thm}
The Hochschild and cyclic homology of $A=A(\mathcal{A}^{(n)\ast},W)$, $n \geq 4$, where $W$ is any cell system on $\mathcal{A}^{(n)}$, is given by
$$\begin{array}{lcl}
HH_0(A) \cong S \oplus C, & \qquad & HC_0(A) \cong S \oplus C, \\
HH_1(A) \cong C, && HC_1(A) = 0, \\
HH_2(A) \cong C^{\ast}[n], && HC_2(A) \cong C^{\ast}[n], \\
HH_3(A) \cong C^{\ast}[n], && HC_3(A) = 0, \\
HH_4(A) \cong C[h], && HC_4(A) \cong C[n], \\
HH_{4+i}(A) \cong HH_i(A)[n], \quad i \geq 1, && HC_{4+i}(A) \cong HC_i(A)[n], \quad i \geq 1,
\end{array}$$
where the graded vector space $C$ has Hilbert series $H_C(t) = \sum_{j=1}^{n-3} \lfloor (n-j-1)/2 \rfloor t^{j}$.
\end{Thm}

\subsubsection{The graph $\mathcal{D}^{(3k)\ast}$}

The Nakayama automorphism is trivial for the graphs $\mathcal{D}^{(3k)\ast}$. We consider the cases $\mathcal{D}^{(6k)\ast}$, $\mathcal{D}^{(6k+3)\ast}$ separately.
For the graph $\mathcal{D}^{(6k)\ast}$, $k \geq 1$, $\mathrm{det}(H_A(t)) = (1-t^6)(1-t^{6k})^{9k-6}/(1-t^3)^{3k-1}$, thus $\chi_{\overline{HC}(A)}(t) = ((3k-1)t^3+(3k-2)t^6+(3k-1)t^9+(3k-2)t^{12}+\cdots+(3k-1)t^{6k-3}-(6k-4)t^{6k})/(1-t^{6k})$. Then since $C$ has Hilbert series $H_C(t) = \sum_{j=1}^{\lfloor (2m-1)/3 \rfloor} (m- \lfloor 3j/2 \rfloor) t^{3j}$, we have $H_X(t) = 0$, $H_K(t) = 6k-4$, and we obtain:

\begin{Thm}
The Hochschild and cyclic homology of $A=A(\mathcal{D}^{(6k)\ast},W)$, $k \geq 1$, where $W$ is equivalent to one of the cell systems constructed in \cite{evans/pugh:2009i}, is given by
$$\begin{array}{lcl}
HH_0(A) \cong S \oplus C, & \qquad & HC_0(A) \cong S \oplus C, \\
HH_1(A) \cong C, && HC_1(A) = 0, \\
HH_2(A) \cong C^{\ast}[6k], && HC_2(A) \cong C^{\ast}[6k], \\
HH_3(A) \cong C^{\ast}[6k] \oplus K[6k], && HC_3(A) \cong K[6k], \\
HH_4(A) \cong C[6k] \oplus K^{\ast}[6k], && HC_4(A) \cong C[6k], \\
HH_{4+i}(A) \cong HH_i(A)[6k], \quad i \geq 1, && HC_{4+i}(A) \cong HC_i(A)[6k], \quad i \geq 1,
\end{array}$$
where the graded vector spaces $C$, $K$ have Hilbert series $H_C(t) = \sum_{j=1}^{\lfloor (2m-1)/3 \rfloor} (m- \lfloor 3j/2 \rfloor) t^{3j}$ and $H_K(t)=6k-4$ respectively.
\end{Thm}

For $\mathcal{D}^{(6k+3)\ast}$, $k \geq 1$, $\mathrm{det}(H_A(t)) = (1-t^{6k+3})^{9k}/(1-t^3)^{3k}$, thus $\chi_{\overline{HC}(A)}(t) = (3kt^3+3kt^6+\cdots+3kt^{6k}-6kt^{6k+3})/(1-t^{6k+3})$. Since $C$ has Hilbert series $H_C(t) = \sum_{j=1}^{\lfloor (2m-2)/3 \rfloor} (m- \lfloor (3j+1)/2 \rfloor) t^{3j}$, we have $H_X(t) = 0$, $H_K(t) = 6k$, and we obtain:

\begin{Thm}
The Hochschild and cyclic homology of $A=A(\mathcal{D}^{(6k+3)\ast},W)$, $k \geq 1$, where $W$ is equivalent to one of the cell systems constructed in \cite{evans/pugh:2009i}, is given by
$$\begin{array}{lcl}
HH_0(A) \cong S \oplus C, & \qquad & HC_0(A) \cong S \oplus C, \\
HH_1(A) \cong C, && HC_1(A) \cong X, \\
HH_2(A) \cong C^{\ast}[6k+3], && HC_2(A) \cong C^{\ast}[6k+3], \\
HH_3(A) \cong C^{\ast}[6k+3] \oplus K[6k+3], && HC_3(A) \cong K[6k+3], \\
HH_4(A) \cong C[6k+3] \oplus K^{\ast}[6k+3], && HC_4(A) \cong C[6k+3], \\
HH_{4+i}(A) \cong HH_i(A)[6k+3], \quad i \geq 1, && HC_{4+i}(A) \cong HC_i(A)[6k+3], \quad i \geq 1,
\end{array}$$
where the graded vector spaces $C$, $K$ have Hilbert series $H_C(t) = \sum_{j=1}^{\lfloor (2m-2)/3 \rfloor} (m- \lfloor (3j+1)/2 \rfloor) t^{3j}$ and $H_K(t)=6k$ respectively.
\end{Thm}

\subsubsection{The graph $\mathcal{E}^{(8)\ast}$}

For the graph $\mathcal{E}^{(8)\ast}$, $\mathrm{det}(H_A(t)) = (1-t^2)(1-t^4)(1-t^8)^2/(1-t)^2$, thus $\chi_{\overline{HC}(A)}(t) = (2t+t^2+2t^3+2t^5+t^6+2t^7-2t^8)/(1-t^8)$. Then $H_X(t) = 0$, $H_K(t) = 2$, and we obtain:

\begin{Thm}
The Hochschild and cyclic homology of $A=A(\mathcal{E}^{(8)\ast},W)$, where $W$ is any cell system on $\mathcal{E}^{(8)\ast}$, is given by
$$\begin{array}{lcl}
HH_0(A) \cong S \oplus C, & \qquad & HC_0(A) \cong S \oplus C, \\
HH_1(A) \cong C, && HC_1(A) = 0, \\
HH_2(A) \cong C^{\ast}[8], && HC_2(A) \cong C^{\ast}[8], \\
HH_3(A) \cong C^{\ast}[8] \oplus K[8], && HC_3(A) \cong K[8], \\
HH_4(A) \cong C[8] \oplus K^{\ast}[8], && HC_4(A) \cong C[8], \\
HH_{4+i}(A) \cong HH_i(A)[8], \quad i \geq 1, && HC_{4+i}(A) \cong HC_i(A)[8], \quad i \geq 1,
\end{array}$$
where the graded vector spaces $C$, $K$ have Hilbert series $H_C(t) = 2t + t^2 + t^3 + t^5$ and $H_K(t)=2$ respectively.
\end{Thm}

\subsection{Determining the Hochschild homology of $A(\mathcal{G},W)$ for non-trivial Nakayama automorphism} \label{Sect:HH(A)-non-trivial_beta}

We now determine the Hochschild and cyclic homology for the graphs $\mathcal{A}^{(n)}$, $n=4,5,6,7$, $\mathcal{E}^{(8)}$.
Here the almost Calabi-Yau algebra $A$ has non-trivial Nakayama automorphism.

By a similar argument to that used in Section \ref{Sect:HH(A)-trivial_beta}, for any almost Calabi-Yau algebra $A$ with non-trivial Nakayama automorphism, we have
$$\begin{array}{cccccccc}
&& 0 &&&& \qquad & \\
(\textrm{min deg},\textrm{max deg}) && \downarrow &&&&& \\
(0,h-3) && \overline{HH}_0(A) & \cong & C &&& \overline{HC}_0(A) \cong C \\
&& \hspace{-3mm} {\scriptstyle B} \downarrow && \hspace{-2mm} {\scriptstyle \cong} \downarrow &&& \\
(1,h-2) && \overline{HH}_1(A) & \cong & C & \hspace{-3mm} \oplus X_1 && \overline{HC}_1(A) \cong X_1 \\
&& \hspace{-3mm} {\scriptstyle B} \downarrow &&& \hspace{-3mm} {\scriptstyle \cong} \downarrow && \\
(2,h-1) && \overline{HH}_2(A) & \cong & X_2 & \hspace{-3mm} \oplus X_1 && \overline{HC}_2(A) \cong X_2 \\
&& \hspace{-3mm} {\scriptstyle B} \downarrow && \hspace{-2mm} {\scriptstyle \cong} \downarrow &&& \\
(3,h) && \overline{HH}_3(A) & \cong & X_2 & \hspace{-3mm} \oplus K_1[h] && \overline{HC}_3(A) \cong K_1[h] \\
&& \hspace{-3mm} {\scriptstyle B} \downarrow &&& \hspace{-3mm} {\scriptstyle \cong} \downarrow && \\
(h,2h-3) && \overline{HH}_4(A) & \cong & X_3 & \hspace{-3mm} \oplus K_1[h] && \overline{HC}_4(A) \cong X_3 \\
&& \hspace{-3mm} {\scriptstyle B} \downarrow && \hspace{-2mm} {\scriptstyle \cong} \downarrow &&& \\
(h+1,2h-2) && \overline{HH}_5(A) & \cong & X_3 & \hspace{-3mm} \oplus X_4 && \overline{HC}_5(A) \cong X_4 \\
&& \hspace{-3mm} {\scriptstyle B} \downarrow &&& \hspace{-3mm} {\scriptstyle \cong} \downarrow && \\
(h+2,2h-1) && \overline{HH}_6(A) & \cong & X_3^{\ast}[3h] & \hspace{-3mm} \oplus X_4^{\ast}[3h] && \overline{HC}_6(A) \cong X_3^{\ast}[3h] \\
&& \hspace{-3mm} {\scriptstyle B} \downarrow && \hspace{-2mm} {\scriptstyle \cong} \downarrow &&& \\
(h+3,2h) && \overline{HH}_7(A) & \cong & X_3^{\ast}[3h] & \hspace{-3mm} \oplus K_1^{\ast}[2h] && \overline{HC}_7(A) \cong K_1^{\ast}[2h] \\
&& \hspace{-3mm} {\scriptstyle B} \downarrow &&& \hspace{-3mm} {\scriptstyle \cong} \downarrow && \\
(2h,3h-3) && \overline{HH}_8(A) & \cong & X_2^{\ast}[3h] & \hspace{-3mm} \oplus K_1^{\ast}[2h] && \overline{HC}_8(A) \cong X_2^{\ast}[3h] \\
&& \hspace{-3mm} {\scriptstyle B} \downarrow && \hspace{-2mm} {\scriptstyle \cong} \downarrow &&& \\
(2h+1,3h-2) && \overline{HH}_{9}(A) & \cong & X_2^{\ast}[3h] & \hspace{-3mm} \oplus X_1^{\ast}[3h] && \overline{HC}_9(A) \cong X_1^{\ast}[3h] \\
&& \hspace{-3mm} {\scriptstyle B} \downarrow &&& \hspace{-3mm} {\scriptstyle \cong} \downarrow && \\
(2h+2,3h-1) && \overline{HH}_{10}(A) & \cong & C^{\ast}[3h] & \hspace{-3mm} \oplus X_1^{\ast}[3h] && \overline{HC}_{10}(A) \cong C^{\ast}[3h] \\
&& \hspace{-3mm} {\scriptstyle B} \downarrow && \hspace{-2mm} {\scriptstyle \cong} \downarrow &&& \\
(2h+3,3h) && \overline{HH}_{11}(A) & \cong & C^{\ast}[3h] & \hspace{-3mm} \oplus K_2[3h] && \overline{HC}_{11}(A) \cong K_2[3h] \\
&& \hspace{-3mm} {\scriptstyle B} \downarrow &&& \hspace{-3mm} {\scriptstyle \cong} \downarrow && \\
(3h,4h-3) && \overline{HH}_{12}(A) & \cong & C[3h] & \hspace{-3mm} \oplus K_2^{\ast}[3h] && \overline{HC}_{12}(A) \cong C[3h] \\
&& \downarrow &&&&& \\
&& \vdots &&&&& \\
\end{array}$$
where $X_1$ lives in degrees 2 to $h-2$, $X_2$ lives in degrees 3 to $h-1$, $X_3$ lives in degrees $h+1$ to $2h-3$, $X_4$ lives in degrees $h+2$ to $2h-2$, $K_i$ lives in degree 0, $i=1,2$, and $\overline{HH}_{13+i}(A) \cong \overline{HH}_{1+i}(A)[3h]$, $\overline{HC}_{12+i}(A) \cong \overline{HC}_i(A)[3h]$ for $i \geq 0$.
The graded vector space $K_1$ can be determined from the Euler characteristic $\chi_{\overline{HC}(A)}(t)$ as it is the only vector space which lives in degree $h$. The vector spaces $X_1$, $X_3$ can be determined by computing $\overline{HH}_1(A)$, $\overline{HH}_4(A)$ respectively. Then $X_2$, $X_4$, $K_2$ can each be determined from knowledge of $C \cong \overline{HH}_0(A)$, $X_1$, $X_3$ and the Euler characteristic.

\subsubsection{The $\mathcal{A}$ graphs}

Here we determine the Hochschild and cyclic homology for the graphs $\mathcal{A}^{(n)}$, $n=4,5,6,7$. The graphs $\mathcal{A}^{(n)}$, $n=5,6,7$, are illustrated in Figure \ref{fig:A(5,6,7)}. We have not yet been able to determine the Hochschild and cyclic homology for the case of general $n$.

\begin{figure}[tb]
\begin{center}
  \includegraphics[width=120mm]{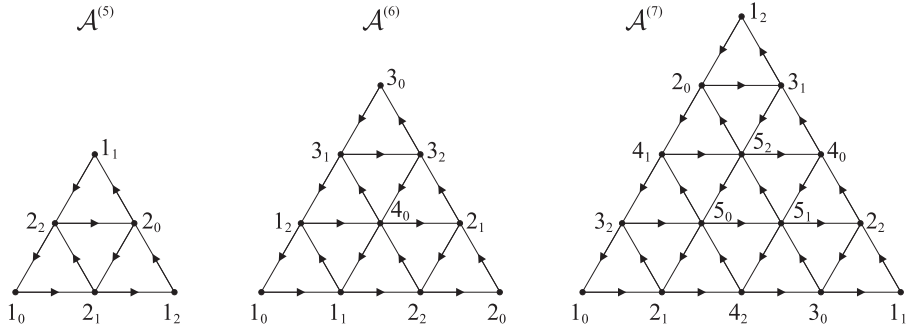}\\
  \caption{Graphs $\mathcal{A}^{(n)}$, $n=5,6,7$} \label{fig:A(5,6,7)}
\end{center}
\end{figure}

We first consider the graph $\mathcal{A}^{(4)}$, for which $\mathrm{det}(H_A(t)) = (1-t^{6})/(1-t^3)$. Thus $\chi_{\overline{HC}(A)}(t) = (t^3+t^6)/(1-t^{12})$ and we see that $H_{K_1}(t) = 0$, and since $C=0$ for all the $\mathcal{A}$ graphs, $H_{K_2}(t) = 0$. Since $\mathrm{ker}(\mu_1') \subset (V \otimes_S A)^S = 0$ and $\mathrm{ker}(\mu_4') = \mathcal{N}^S = 0$, we see that $\overline{HH}_1(A) = 0 = \overline{HH}_4(A)$. Thus $X_1 = X_3 = K_1 = 0$, and from $\chi_{\overline{HC}(A)}(t)$ we deduce that $X_2$ has Hilbert series $H_{X_2}(t) = t^3$ and $X_4 = K_2 = 0$.

\begin{Thm}
The Hochschild and cyclic homology of $A=A(\mathcal{A}^{(4)},W)$, where $W$ is any cell system on $\mathcal{A}^{(4)}$, is given by
$$\begin{array}{lcl}
HH_0(A) \cong S, & \qquad & HC_0(A) \cong S, \\
HH_1(A) = 0, && HC_1(A) = 0, \\
HH_2(A) \cong X, && HC_2(A) \cong X, \\
HH_3(A) \cong X, && HC_3(A) = 0, \\
HH_4(A) = 0, && HC_4(A) = 0, \\
HH_5(A) = 0, && HC_5(A) = 0, \\
HH_6(A) = 0, && HC_6(A) = 0, \\
HH_7(A) = 0, && HC_7(A) = 0, \\
HH_8(A) \cong X^{\ast}[12], && HC_8(A) \cong X^{\ast}[12], \\
HH_9(A) \cong X^{\ast}[12], && HC_9(A) = 0, \\
HH_{10}(A) = 0, && HC_{10}(A) = 0, \\
HH_{11}(A) = 0, && HC_{11}(A) = 0, \\
HH_{12}(A) = 0, && HC_{12}(A) = 0, \\
HH_{12+i}(A) \cong HH_i(A)[12], \quad i \geq 1, && HC_{12+i}(A) \cong HC_i(A)[12], \quad i \geq 1,
\end{array}$$
where the graded vector space $X$ has Hilbert series $H_{X}(t) = t^3$.
\end{Thm}

We now consider the graph $\mathcal{A}^{(5)}$, for which $\mathrm{det}(H_A(t)) = (1-t^{15})/(1-t^3)$. Thus $\chi_{\overline{HC}(A)}(t) = (t^3+t^6+t^9+t^{12})/(1-t^{15})$ and we see that $H_{K_1}(t) = 0$, and since $C=0$ for all the $\mathcal{A}$ graphs, $H_{K_2}(t) = 0$. We now explicitly determine $\overline{HH}_1(A)$ and $\overline{HH}_4(A)$.

We begin with the graded vector space $Y = \overline{HH}_1(A) = \mathrm{ker}(\mu_1')/\mathrm{Im}(\mu_2')$, and consider each graded piece $Y_j$ separately. Due to the three-colourability of $\mathcal{A}^{(5)}$, $Y_j = 0$ for $j=1,2$. Thus we only need to determine $Y_3$. A basis for $(\widetilde{V} \otimes_S A)^S_2$ is given by the elements $[2_{l+1} 1_l] \otimes [1_l 2_{l+1}]$, $[1_{l+1} 2_l] \otimes [2_l 1_{l+1}]$ and $[2_{l+1} 2_l] \otimes [2_l 2_{l+1}]$, for $l=0,1,2 \textrm{ mod } 3$.
We have $\mu_2'([2_{l+1} 1_l] \otimes [1_l 2_{l+1}]) = W_{1_l 2_{l+1} 2_{l+2}} ([2_{l+1} 2_{l+2}] \otimes [2_{l+2} 1_l 2_{l+1}] + [1_l 2_{l+1}] \otimes [2_{l+1} 2_{l+2} 1_l]) = W_{2_l 2_{l+1} 2_{l+2}} [2_{l+1} 2_{l+2}] \otimes [2_{l+2} 2_l 2_{l+1}]$, using the relations in $A$.
Thus we see that $[2_{l+1} 2_{l+2}] \otimes [2_{l+2} 2_l 2_{l+1}] = 0$ in $Y_3$, $l=0,1,2 \textrm{ mod } 3$. Thus $Y_3 = 0$ and we obtain $\overline{HH}_1(A) = 0$. Since $X_1 = 0$, we deduce from $\chi_{\overline{HC}(A)}(t)$ that $X_2$ has Hilbert series $H_{X_2}(t) = t^3$.

We now consider $Y' = \overline{HH}_4(A)$, which lives in degrees 5 to 7. Now $\mathrm{ker}(\mu_4') = \mathcal{N}^S[5]$, since $\sum_j w_j^{\ast} \beta(x) w_j = 0$ for all $x \in \mathcal{N}^S_+$ and $\mathcal{N}^S_0 = 0$. As with $\overline{HH}_1(A)$, $Y'_j = 0$ for $j=5,7$, due to the three-colourability of $\mathcal{A}^{(5)}$.
We now determine $Y'_6$. A basis for $(V \otimes_S \mathcal{N})^S_1$ is given by $[2_l 2_{l+1}] \otimes [2_{l+1}]$, $l=0,1,2 \textrm{ mod } 3$, and a basis for $\mathcal{N}^S$ is given by $[2_l 2_{l+1}]$, $l=0,1,2 \textrm{ mod } 3$. Now $\mu_5'([2_l 2_{l+1}] \otimes [2_{l+1}]) = [2_{l+1} 2_{l+2}] - [2_l 2_{l+1}]$, thus $[2_0 2_1] = [2_1 2_2] = [2_2 2_0]$ in $Y'_6$. Then $\overline{HH}_4(A) = \mathbb{C}[2_0 2_1] [5] = X_3$, and we deduce from $\chi_{\overline{HC}(A)}(t)$ that $X_4 = 0$.

\begin{Thm}
The Hochschild and cyclic homology of $A=A(\mathcal{A}^{(5)},W)$, where $W$ is any cell system on $\mathcal{A}^{(5)}$, is given by
$$\begin{array}{lcl}
HH_0(A) \cong S, & \qquad & HC_0(A) \cong S, \\
HH_1(A) = 0, && HC_1(A) = 0, \\
HH_2(A) \cong X_2, && HC_2(A) \cong X_2, \\
HH_3(A) \cong X_2, && HC_3(A) = 0, \\
HH_4(A) \cong X_3, && HC_4(A) \cong X_3, \\
HH_5(A) \cong X_3, && HC_5(A) = 0, \\
HH_6(A) \cong X_3^{\ast}[15], && HC_6(A) \cong X_3^{\ast}[15], \\
HH_7(A) \cong X_3^{\ast}[15], && HC_7(A) = 0, \\
HH_8(A) \cong X_2^{\ast}[15], && HC_8(A) \cong X_2^{\ast}[15], \\
HH_9(A) \cong X_2^{\ast}[15], && HC_9(A) = 0, \\
HH_{10}(A) = 0, && HC_{10}(A) = 0, \\
HH_{11}(A) = 0, && HC_{11}(A) = 0, \\
HH_{12}(A) = 0, && HC_{12}(A) = 0, \\
HH_{12+i}(A) \cong HH_i(A)[15], \quad i \geq 1, && HC_{12+i}(A) \cong HC_i(A)[15], \quad i \geq 1,
\end{array}$$
where the graded vector spaces $X_2$ and $X_3$ have Hilbert series $H_{X_2}(t) = t^3$ and $H_{X_3}(t) = t^6$.
\end{Thm}

We now consider the graph $\mathcal{A}^{(6)}$, for which $\mathrm{det}(H_A(t)) = (1-t^6)(1-t^9)(1-t^{18})/(1-t^3)$. Thus $\chi_{\overline{HC}(A)}(t) = (t^3+t^{15}-2t^{18})/(1-t^{18})$ and we see that $H_{K_1}(t) = 0$, and since $C=0$ for all the $\mathcal{A}$ graphs, $H_{K_2}(t) = 2$. We now explicitly determine $\overline{HH}_1(A)$ and $\overline{HH}_4(A)$.

We begin with the graded vector space $Y = \overline{HH}_1(A)$. Due to the three-colourability of $\mathcal{A}^{(6)}$, $Y_j = 0$ for $j=1,2,4$, so we only need to determine $Y_3$. A basis for $(\widetilde{V} \otimes_S A)^S_2$ is given by the elements $[i_{l+1} i_l] \otimes [i_l i_{l+1}]$, $[(i+1)_2 i_1] \otimes [i_1 (i+1)_2]$, $[i_1 4_0] \otimes [4_0 i_1]$, and $[4_0 i_2] \otimes [i_2 4_0]$, for $l=0,1,2$, $i=1,2,3 \textrm{ mod } 3$.
A basis for $(V \otimes A)^S_3$ is given by $[4_0 i_1] \otimes [i_1 i_2 4_0]$, $[i_1 i_2] \otimes [i_2 4_0 i_1]$ and $[i_2 4_0] \otimes [4_0 i_1 i_2]$, for $i=1,2,3 \textrm{ mod } 3$.
Under $\mu_2'$ the basis elements of $(\widetilde{V} \otimes_S A)^S_2$ yield the following expressions after using the relations in $A$, for $i=1,2,3 \textrm{ mod } 3$:
\begin{eqnarray*}
&& \mu_2'([i_1 i_0] \otimes [i_0 i_1]) = -W_{4_0 i_1 i_2} [i_1 i_2] \otimes [i_2 4_0 i_1] = \mu_2'([i_0 i_2] \otimes [i_2 i_0]), \\
&& \mu_2'([i_2 i_1] \otimes [i_1 i_2]) = W_{4_0 i_1 i_2} [4_0 i_1] \otimes [i_1 i_2 4_0] + W_{4_0 i_1 i_2} [i_2 4_0] \otimes [4_0 i_1 i_2], \\
&& \mu_2'([(i+1)_2 i_1] \otimes [i_1 (i+1)_2]) = W_{4_0 i_1 (i+1)_2} [4_0 i_1] \otimes [i_1 (i+1)_2 4_0] \\
&& \qquad \qquad \qquad \qquad \qquad \qquad \qquad + W_{4_0 i_1 (i+1)_2} [(i+1)_2 4_0] \otimes [4_0 i_1 (i+1)_2], \\
&& \mu_2'([i_1 4_0] \otimes [4_0 i_1]) = W_{4_0 i_1 i_2} [i_2 4_0] \otimes [4_0 i_1 i_2] + W_{4_0 i_1 (i+1)_2} [(i+1)_2 4_0] \otimes [4_0 i_1 (i+1)_2] \\
&& \qquad \qquad \qquad \qquad \qquad + W_{4_0 i_1 i_2} [i_1 i_2] \otimes [i_2 4_0 i_1], \\
&& \mu_2'([4_0 i_2] \otimes [i_2 4_0]) = W_{4_0 i_1 i_2} [4_0 i_1] \otimes [i_1 i_2 4_0] + W_{4_0 (i-1)_1 i_2} [4_0 (i-1)_1] \otimes [(i-1)_1 i_2 4_0] \\
&& \qquad \qquad \qquad \qquad \qquad + W_{4_0 i_1 i_2} [i_1 i_2] \otimes [i_2 4_0 i_1].
\end{eqnarray*}
Then from $\mathrm{Im}(\mu_2')_3$ we obtain the following relations in $Y_3$: $[i_1 i_2] \otimes [i_2 4_0 i_1] = 0$ and $[4_0 i_1] \otimes [i_1 i_2 4_0] = -[i_2 4_0] \otimes [4_0 i_1 i_2] = (W_{4_0 3_1 3_2}/W_{4_0 i_1 i_2}) [4_0 3_1] \otimes [3_1 1_2 4_0]$.

We now consider $\mathrm{Ker}(\mu_1')_3$. Let $x = \sum_{i=1}^3 (\lambda_i^0 [i_2 4_0] \otimes [4_0 i_1 i_2] + \lambda_i^1 [4_0 i_1] \otimes [i_1 i_2 4_0] + \lambda_i^2 [i_1 i_2] \otimes [i_2 4_0 i_1])$ be a general element in $(V \otimes A)^S_3$. Since $\mu_1'(x) = \sum_{i=1}^3((\lambda_i^0 - \lambda_i^1) [4_0 i_1 i_2 4_0] + (\lambda_i^1 - \lambda_i^2) [i_1 i_2 4_0 i_1] + (\lambda_i^2 - \lambda_i^0) [i_2 4_0 i_1 i_2])$, then $x \in \mathrm{Ker}(\mu_1')$ if and only if $\lambda_i^0 = \lambda_i^1 = \lambda_i^2$ for each $i=1,2,3$. Using the relations from $\mathrm{Im}(\mu_2')_3$, a general element in $Y_3$ is thus of the form $\sum_i \lambda_i^0 ([i_2 4_0] \otimes [4_0 i_1 i_2] + [4_0 i_1] \otimes [i_1 i_2 4_0] + [i_1 i_2] \otimes [i_2 4_0 i_1]) = \sum_i \lambda_i^0 (-(W_{4_0 3_1 3_2}/W_{4_0 i_1 i_2}) [4_0 3_1] \otimes [3_1 1_2 4_0] + 0 + (W_{4_0 3_1 3_2}/W_{4_0 i_1 i_2}) [4_0 3_1] \otimes [3_1 1_2 4_0]) = 0$. Thus $Y_3 = 0$ and we obtain $\overline{HH}_1(A) = 0$. Since $X_1 = 0$, we deduce from $\chi_{\overline{HC}(A)}(t)$ that $X_2$ is a graded vector space with Hilbert series $H_{X_2}(t) = t^3$.

We now consider $Y' = \overline{HH}_4(A)$, which lives in degrees 6 to 9.
Now $\mathcal{N}^S_0 = \mathbb{C}[4_0]$, but $\mu_4'(\lambda [4_0]) = 2\lambda [4_0 1_1 1_2 4_0]$, thus $\lambda [4_0] \in \mathrm{ker}(\mu_4')$ if and only if $\lambda = 0$.
Thus $\mathrm{ker}(\mu_4') = \mathcal{N}^S_+[6]$, since $\sum_j w_j^{\ast} \beta(x) w_j = 0$ for all $x \in \mathcal{N}^S_+ = \mathbb{C} [4_0 1_1 1_2 4_0]$.
As with $\overline{HH}_1(A)$, $Y'_j = 0$ for $j=7,8$, due to the three-colourability of $\mathcal{A}^{(6)}$, and $Y'_6 = 0$ since $\mathrm{ker}(\mu_4')_6 = 0$.
We now determine $Y'_9$. A basis for $(V \otimes_S \mathcal{N})^S_1$ is given by $[4_0 i_1] \otimes [i_1 i_2 4_0]$, $[i_1 (i+1)_2] \otimes [(i+1)_2 4_0 (i-1)_1]$ and $[i_2 4_0] \otimes [4_0 (i-1)_1 (i-1)_2]$, $i=0,1,2 \textrm{ mod } 3$.
Now $\mu_5'([4_0 1_1] \otimes [1_1 1_2 4_0]) = -[4_0 1_1 1_2 4_0]$, thus $\mathrm{Im}(\mu_5') = \mathcal{N}^S_+$. Then $Y'_9 = 0$ and we obtain $\overline{HH}_4(A) = 0$.
Then since $X_3 = 0$, we deduce from $\chi_{\overline{HC}(A)}(t)$ that $X_4 = 0$.

\begin{Thm}
The Hochschild and cyclic homology of $A=A(\mathcal{A}^{(6)},W)$, where $W$ is any cell system on $\mathcal{A}^{(6)}$, is given by
$$\begin{array}{lcl}
HH_0(A) \cong S, & \qquad & HC_0(A) \cong S, \\
HH_1(A) = 0, && HC_1(A) = 0, \\
HH_2(A) \cong X, && HC_2(A) \cong X, \\
HH_3(A) \cong X, && HC_3(A) = 0, \\
HH_4(A) = 0, && HC_4(A) = 0, \\
HH_5(A) = 0, && HC_5(A) = 0, \\
HH_6(A) = 0, && HC_6(A) = 0, \\
HH_7(A) = 0, && HC_7(A) = 0, \\
HH_8(A) \cong X^{\ast}[18], && HC_8(A) \cong X^{\ast}[18], \\
HH_9(A) \cong X^{\ast}[18], && HC_9(A) = 0, \\
HH_{10}(A) = 0, && HC_{10}(A) = 0, \\
HH_{11}(A) \cong K[18], && HC_{11}(A) \cong K[18], \\
HH_{12}(A) \cong K^{\ast}[18], && HC_{12}(A) = 0, \\
HH_{12+i}(A) \cong HH_i(A)[18], \quad i \geq 1, && HC_{12+i}(A) \cong HC_i(A)[18], \quad i \geq 1,
\end{array}$$
where the graded vector spaces $X$ and $K$ have Hilbert series $H_{X}(t) = t^3$ and $H_K(t) = 2$.
\end{Thm}

We now consider the graph $\mathcal{A}^{(7)}$, for which $\mathrm{det}(H_A(t)) = (1-t^{21})^3/(1-t^3)$. Thus $\chi_{\overline{HC}(A)}(t) = (t^3+t^6+\cdots+t^{18}-2t^{21})/(1-t^{21})$ and we see that $H_{K_1}(t) = 0$, and since $C=0$ for all the $\mathcal{A}$ graphs, $H_{K_2}(t) = 2$. We now explicitly determine $\overline{HH}_1(A)$ and $\overline{HH}_4(A)$.

We begin with the graded vector space $Y = \overline{HH}_1(A)$. Due to the three-colourability of $\mathcal{A}^{(7)}$, $Y_j = 0$ for $j=1,2,4,5$, so we only need to determine $Y_3$. A basis for $(\widetilde{V} \otimes_S A)^S_2$ is given by the elements $[1_l 3_{l-1}] \otimes [3_{l-1} 1_l]$, $[2_l 1_{l-1}] \otimes [1_{l-1} 2_l]$, $[2_l 5_{l-1}] \otimes [5_{l-1} 2_l]$, $[3_l 2_{l-1}] \otimes [2_{l-1} 3_l]$, $[3_l 4_{l-1}] \otimes [4_{l-1} 3_l]$, $[4_l 2_{l-1}] \otimes [2_{l-1} 4_l]$, $[4_l 5_{l-1}] \otimes [5_{l-1} 4_l]$, $[5_l 3_{l-1}] \otimes [3_{l-1} 5_l]$, $[5_l 4_{l-1}] \otimes [4_{l-1} 5_l]$ and $[5_l 5_{l-1}] \otimes [5_{l-1} 5_l]$, for $l=0,1,2 \textrm{ mod } 3$.
A basis for $(V \otimes A)^S_3$ is given by $[2_l 3_{l+1}] \otimes [3_{l+1} 5_{l+2} 2_l]$, $[3_l 5_{l+1}] \otimes [5_{l+1} 2_{l+2} 3_l]$, $[4_l 5_{l+1}] \otimes [5_{l+1} 2_{l+2} 4_l]$, $[5_l 2_{l+1}] \otimes [2_{l+1} 3_{l+2} 5_l]$, $[5_l 4_{l+1}] \otimes [4_{l+1} 3_{l+2} 5_l]$ and $[5_l 5_{l+1}] \otimes [5_{l+1} 4_{l+2} 5_l]$, for $l=0,1,2 \textrm{ mod } 3$.
Under $\mu_2'$ the basis elements of $(\widetilde{V} \otimes_S A)^S_2$ yield the following expressions after using the relations in $A$, for $l=0,1,2 \textrm{ mod } 3$:
\begin{eqnarray*}
\mu_2'([1_l 3_{l-1}] \otimes [3_{l-1} 1_l]) & = & W_{235} [2_{l+1} 3_{l-1}] \otimes [3_{l-1} 5_l 2_{l+1}] \;\; = \;\; \mu_2'([2_{l+1} 1_l] \otimes [1_l 2_{l+1}]), \\
\mu_2'([2_l 5_{l-1}] \otimes [5_{l-1} 2_l]) & = & W_{235} [3_{l+1} 5_{l-1}] \otimes [5_{l-1} 2_l 3_{l+1}] + W_{245} [4_{l+1} 5_{l-1}] \otimes [5_{l-1} 2_l 4_{l+1}] \\
&& \quad + W_{235} [2_{l+1} 3_{l-1}] \otimes [3_{l-1} 5_l 2_{l+1}], \\
\mu_2'([3_l 2_{l-1}] \otimes [2_{l-1} 3_l]) & = & W_{235} [5_{l+1} 2_{l-1}] \otimes [2_{l-1} 3_l 5_{l+1}] + W_{235} [3_l 5_{l+1}] \otimes [5_{l+1} 2_{l-1} 3_l], \\
\mu_2'([3_l 4_{l-1}] \otimes [4_{l-1} 3_l]) & = & W_{354} [5_{l+1} 4_{l-1}] \otimes [4_{l-1} 3_l 5_{l+1}] - W_{235} [3_l 5_{l+1}] \otimes [5_{l+1} 2_{l-1} 3_l], \\
\mu_2'([4_l 2_{l-1}] \otimes [2_{l-1} 4_l]) & = & -W_{235} [5_{l+1} 2_{l-1}] \otimes [2_{l-1} 3_l 5_{l+1}] + W_{245} [4_l 5_{l+1}] \otimes [5_{l+1} 2_{l-1} 4_l], \\
\mu_2'([4_l 5_{l-1}] \otimes [5_{l-1} 4_l]) & = & W_{455} [5_{l+1} 5_{l-1}] \otimes [5_{l-1} 4_l 5_{l+1}] - W_{245} [4_l 5_{l+1}] \otimes [5_{l+1} 2_{l-1} 4_l] \\
&& \quad - W_{235} [3_{l+1} 5_{l-1}] \otimes [5_{l-1} 2_l 3_{l+1}], \\
\mu_2'([5_l 3_{l-1}] \otimes [3_{l-1} 5_l]) & = & W_{354} [5_l 4_{l+1}] \otimes [4_{l+1} 3_{l-1} 5_l] + W_{235} [5_l 2_{l+1}] \otimes [2_{l+1} 3_{l-1} 5_l] \\
&& \quad + W_{235} [2_{l+1} 3_{l-1}] \otimes [3_{l-1} 5_l 2_{l+1}], \\
\mu_2'([5_l 4_{l-1}] \otimes [4_{l-1} 5_l]) & = & -W_{354} [5_{l+1} 4_{l-1}] \otimes [4_{l-1} 3_l 5_{l+1}] + W_{455} [5_l 5_{l+1}] \otimes [5_{l+1} 4_{l-1} 5_l] \\
&& \quad - W_{235} [5_l 2_{l+1}] \otimes [2_l 3_{l-1} 5_l], \\
\mu_2'([5_l 5_{l-1}] \otimes [5_{l-1} 5_l]) & = & -W_{245} [4_{l+1} 5_{l-1}] \otimes [5_{l-1} 2_l 4_{l+1}] - W_{354} [5_l 4_{l+1}] \otimes [4_{l+1} 3_{l-1} 5_l] \\
&& \quad - W_{455} [5_{l+1} 5_{l-1}] \otimes [5_{l-1} 4_l 5_{l+1}] - W_{455} [5_l 5_{l+1}] \otimes [5_{l+1} 4_{l-1} 5_l],
\end{eqnarray*}
where $W_{ijk} := W_{i_0 j_1 k_2} = W_{i_1 j_2 k_0} = W_{i_2 j_0 k_1}$ for vertices $i_l$, $j_l$, $k_l$ of $\mathcal{A}^{(7)}$, $l=0,1,2$.
Then from $\mathrm{Im}(\mu_2')_3$ we obtain the following relations in $Y_3$:
$$\begin{array}{l}
{} [2_{l} 3_{l+1}] \otimes [3_{l+1} 5_{l+2} 2_{l}] = 0, \quad l=0,1,2 \textrm{ mod } 3, \qquad \qquad [5_0 5_1] \otimes [5_1 4_2 5_0] = -x_2 - x_3, \\
{} [3_0 5_1] \otimes [5_1 2_2 3_0] = -(W_{245}/W_{235}) [4_0 5_1] \otimes [5_1 2_2 4_0] = -[5_1 2_2] \otimes [2_2 3_0 5_1] \\
\hspace{29.5mm} = (W_{354}/W_{235}) [5_1 4_2] \otimes [4_2 3_0 5_1] = (W_{354}/W_{235}) x_1 - (W_{455}/W_{235}) x_2, \\
{} [3_1 5_2] \otimes [5_2 2_0 3_1] = -(W_{245}/W_{235}) [4_1 5_2] \otimes [5_2 2_0 4_1] = -[5_2 2_0] \otimes [2_0 3_1 5_2] = (W_{354}/W_{235}) x_1, \\
{} [3_2 5_0] \otimes [5_0 2_1 3_2] = -(W_{245}/W_{235}) [4_2 5_0] \otimes [5_0 2_1 4_2] = -[5_0 2_1] \otimes [2_1 3_2 5_0] \\
\hspace{29.5mm} = (W_{354}/W_{235}) [5_0 4_1] \otimes [4_1 3_2 5_0] = (W_{354}/W_{235}) x_1 + (W_{455}/W_{235}) x_3,
\end{array}$$
where $x_1 = [5_2 4_0] \otimes [4_0 3_1 5_2]$, $x_2 = [5_1 5_2] \otimes [5_2 4_0 5_1]$ and $x_3 = [5_2 5_0] \otimes [5_0 4_1 5_2]$.

We now consider $\mathrm{Ker}(\mu_1')_3$. Let $x = \sum_{l=0}^2 (\lambda_1^l [3_l 5_{l+1}] \otimes [5_{l+1} 2_{l+2} 3_l] + \lambda_2^l [4_l 5_{l+1}] \otimes [5_{l+1} 2_{l+2} 4_l] + \lambda_3^l [5_l 2_{l+1}] \otimes [2_{l+1} 3_{l+2} 5_l] + \lambda_4^l [5_l 4_{l+1}] \otimes [4_{l+1} 3_{l+2} 5_l] + \lambda_5^l [5_l 5_{l+1}] \otimes [5_{l+1} 4_{l+2} 5_l])$ be a general element in $(V \otimes A)^S_3$.
Now $\mu_1'(x) = \sum_{l=0}^2 (\lambda_3^{l-1} [2_l 3_{l+1} 5_{l+2} 2_l] - \lambda_1^l [3_l 5_{l+1} 2_{l+2} 3_l] + (\lambda_1^{l-1} - \lambda_2^{l-1} (W_{235}/W_{245}) - \lambda_3^l + \lambda_4^l (W_{235}/W_{354}) + (\lambda_5^{l-1}-\lambda_5^l) (W_{235}/W_{455})) [5_l 2_{l+1} 3_{l+2} 5_l]) = 0$ if and only if $\lambda_1^l = \lambda_3^l = 0, \lambda_4^l = \lambda_2^{l-1} + (\lambda_5^l-\lambda_5^{l-1})(W_{354}/W_{455})$, for $l=0,1,2 \textrm{ mod } 3$.
Using the relations from $\mathrm{Im}(\mu_2')_3$, a general element in $Y_3$ is thus of the form $\sum_{l=0}^2 (\lambda_2^l [4_l 5_{l+1}] \otimes [5_{l+1} 2_{l+2} 4_l] + (\lambda_2^{l-1} + (\lambda_5^l-\lambda_5^{l-1})(W_{354}/W_{455})) [5_l 4_{l+1}] \otimes [4_{l+1} 3_{l+2} 5_l] + \lambda_5^l [5_l 5_{l+1}] \otimes [5_{l+1} 4_{l+2} 5_l]) = 0$.
Thus $Y_3 = 0$ and we obtain $\overline{HH}_1(A) = 0$. Since $X_1 = 0$, we deduce from $\chi_{\overline{HC}(A)}(t)$ that $X_2$ is a graded vector space with Hilbert series $H_{X_2}(t) = t^3 + t^6$.

We now consider $Y' = \overline{HH}_4(A)$, which lives in degrees 7 to 11. Now $\mathrm{ker}(\mu_4') = \mathcal{N}^S[7]$, since $\sum_j w_j^{\ast} \beta(x) w_j = 0$ for all $x \in \mathcal{N}^S_+$ and $\mathcal{N}^S_0 = 0$. As with $\overline{HH}_1(A)$, $Y'_j = 0$ for $j=7,8,10,11$, due to the three-colourability of $\mathcal{A}^{(7)}$.
We now determine $Y'_9$. A basis for $(V \otimes_S \mathcal{N})^S_2$ is given by $[4_l 5_{l+1}] \otimes [5_{l+1} 4_l]$, $[5_l 4_{l+1}] \otimes [4_{l+1} 5_l]$ and $[5_l 5_{l+1}] \otimes [5_{l+1} 5_l]$, and a basis for $\mathcal{N}^S$ is given by $[4_l 5_{l+1} 4_{l+2}]$, $[5_l 5_{l+1} 5_{l+2}]$, for $l=0,1,2 \textrm{ mod } 3$. Using the relations in $A$ we obtain $\mu_5'([4_l 5_{l+1}] \otimes [5_{l+1} 4_l]) = -(W_{555}/W_{455}) [5_{l+1} 5_{l+2} 5_l] - [4_l 5_{l+1} 4_{l+2}]$, $\mu_5'([5_l 4_{l+1}] \otimes [4_{l+1} 5_l]) = (W_{555}/W_{455}) [5_l 5_{l+1} 5_{l+2}] + [4_{l+1} 5_{l+2} 4_l]$ and $\mu_5'([5_l 5_{l+1}] \otimes [5_{l+1} 5_l]) = [5_{l+1} 5_{l+2} 5_l] - [5_l 5_{l+1} 5_{l+2}]$, for $l=0,1,2 \textrm{ mod } 3$. These yield the relations $[4_l 5_{l+1} 4_{l+2}] = -(W_{555}/W_{455}) [5_{l'} 5_{l'+1} 5_{l'+2}]$ in $Y'_9$, for all $l,l'=0,1,2 \textrm{ mod } 3$.
Thus we obtain $\overline{HH}_4(A) = \mathbb{C}[5_0 5_1 5_2] [7]$.
Then $X_3 = \mathbb{C}[5_0 5_1 5_2] [7]$ and we deduce from $\chi_{\overline{HC}(A)}(t)$ that $X_4 = 0$.

\begin{Thm}
The Hochschild and cyclic homology of $A=A(\mathcal{A}^{(7)},W)$, where $W$ is any cell system on $\mathcal{A}^{(7)}$, is given by
$$\begin{array}{lcl}
HH_0(A) \cong S, & \qquad & HC_0(A) \cong S, \\
HH_1(A) = 0, && HC_1(A) = 0, \\
HH_2(A) \cong X_2, && HC_2(A) \cong X_2, \\
HH_3(A) \cong X_2, && HC_3(A) = 0, \\
HH_4(A) \cong X_3, && HC_4(A) \cong X_3, \\
HH_5(A) \cong X_3, && HC_5(A) = 0, \\
HH_6(A) \cong X_3^{\ast}[21], && HC_6(A) \cong X_3^{\ast}[21], \\
HH_7(A) \cong X_3^{\ast}[21], && HC_7(A) = 0, \\
HH_8(A) \cong X_2^{\ast}[21], && HC_8(A) \cong X_2^{\ast}[21], \\
HH_9(A) \cong X_2^{\ast}[21], && HC_9(A) = 0, \\
HH_{10}(A) = 0, && HC_{10}(A) = 0, \\
HH_{11}(A) \cong K[21], && HC_{11}(A) \cong K[21], \\
HH_{12}(A) \cong K^{\ast}[21], && HC_{12}(A) = 0, \\
HH_{12+i}(A) \cong HH_i(A)[21], \quad i \geq 1, && HC_{12+i}(A) \cong HC_i(A)[21], \quad i \geq 1,
\end{array}$$
where the graded vector spaces $X_2$, $X_3$ and $K$ have Hilbert series $H_{X_2}(t) = t^3 + t^6$, $H_{X_3}(t) = t^9$ and $H_K(t) = 2$.
\end{Thm}

\subsubsection{The graph $\mathcal{E}^{(8)}$}

For the graph $\mathcal{E}^{(8)}$, $\mathrm{det}(H_A(t)) = (1-t^6)(1-t^{12})(1-t^{24})^2/(1-t^3)^2 = \mathrm{det}(H_{A'}(t^3))$, where $A' = A(\mathcal{E}^{(8)\ast},W)$. Thus $\chi_{\overline{HC}(A)}(t) = (2t^3+t^6+2t^9+2t^{15}+t^{18}+2t^{21}-2t^{24})/(1-t^{24})$ and we see that $H_{K_1}(t) = 0$. Since $C \cong \mathbb{C} [2_0 2_1 2_2 2_0]$ which lives in degree $>0$, we see that $H_{K_2}(t) = 2$. We now explicitly determine $\overline{HH}_1(A)$ and $\overline{HH}_4(A)$.

We begin with the graded vector space $Y = \overline{HH}_1(A) = \mathrm{ker}(\mu_1')/\mathrm{Im}(\mu_2')$, and consider each graded piece $Y_j$ separately. Due to the three-colourability of $\mathcal{E}^{(8)}$, $Y_j = 0$ for $j=1,2,4,5$. We will first determine $Y_3$. A basis for $(\widetilde{V} \otimes_S A)^S_3$ is given by $[2_{l+1} 1_l] \otimes [1_l 2_{l+1}]$, $[2_{l+1} 2_l] \otimes [2_l 2_{l+1}]$, $[3_{l+1} 2_l] \otimes [2_l 3_{l+1}]$, $[4_{l+1} 2_l] \otimes [2_l 4_{l+1}]$, $[1_{l+1} 3_l] \otimes [3_l 1_{l+1}]$, $[2_{l+1} 3_l] \otimes [3_l 2_{l+1}]$, $[3_{l+1} 3_l] \otimes [3_l 3_{l+1}]$, and $[3_{l+1} 4_l] \otimes [4_l 3_{l+1}]$, for $l=1,2,3$. A basis for $(V \otimes_S A)^S_3$ is given by $[2_{l-1} 2_l] \otimes [2_l 2_{l+1} 2_{l-1}]$, $[3_{l-1} 3_l] \otimes [3_l 3_{l+1} 3_{l-1}]$, $[3_{l-1} 2_l] \otimes [2_l 2_{l+1} 3_{l-1}]$, $[2_{l-1} 3_l] \otimes [3_l 3_{l+1} 2_{l-1}]$, $[3_{l-1} 2_l] \otimes [2_l 3_{l+1} 3_{l-1}]$ and $[2_{l-1} 3_l] \otimes [3_l 2_{l+1} 2_{l-1}]$, $l=1,2,3$.
Under $\mu_2'$, $[2_{l+1} 1_l] \otimes [1_l 2_{l+1}]$ gives
\begin{eqnarray*}
\mu_2'([2_{l+1} 1_l] \otimes [1_l 2_{l+1}]) & = & \sqrt{[2][3]} ([2_{l+1}3_{l-1}] \otimes [3_{l-1}1_l2_{l+1}] + [1_l2_{l+1}] \otimes [2_{l+1}3_{l-1}1_l]) \\
& = & -\sqrt{[3][4]} ([2_{l+1}3_{l-1}] \otimes [3_{l-1}3_l2_{l+1}] + [2_{l+1}3_{l-1}] \otimes [3_{l-1}2_l2_{l+1}]),
\end{eqnarray*}
using the relations in $A$. We get the same result from considering $\mu_2'([1_l 3_{l-1}] \otimes [3_{l-1} 1_l])$. We also obtain (up to some scalar factor)
\begin{eqnarray*}
\mu_2'([3_{l-1} 2_{l+1}] \otimes [2_{l+1} 3_{l-1}]) & = & [3_l2_{l+1}] \otimes [2_{l+1}3_{l-1}3_l] + [3_{l-1}3_l] \otimes [3_l2_{l+1}3_{l-1}] \\
&& + [2_l2_{l+1}] \otimes [2_{l+1}3_{l-1}2_l] + [3_{l-1}2_l] \otimes [2_l2_{l+1}3_{l-1}],
\end{eqnarray*}
and the results for $\mu_2'([3_{l+1} 4_l] \otimes [4_l 3_{l+1}])$, $\mu_2'([4_l 2_{l-1}] \otimes [2_{l-1} 4_l])$ and $\mu_2'([2_{l-1} 3_{l+1}] \otimes [3_{l+1} 2_{l-1}])$ are given by the above results by interchanging $1_p \leftrightarrow 4_p$, $2_p \leftrightarrow 3_p$ for $p=l,l+1,l-1$. Finally, we also have (again up to some scalar factor)
\begin{eqnarray*}
\mu_2'([2_{l-1} 2_{l+1}] \otimes [2_{l+1} 2_{l-1}]) & = & [3_l2_{l+1}] \otimes [2_{l+1}2_{l-1}3_l] + [2_{l-1}3_l] \otimes [3_l2_{l+1}2_{l-1}] \\
&& + \sqrt{[3]}[2_l2_{l+1}] \otimes [2_{l+1}2_{l-1}2_l] + \sqrt{[3]}[2_{l-1}3_l] \otimes [3_l2_{l+1}2_{l-1}], \\
\mu_2'([3_{l-1} 3_{l+1}] \otimes [3_{l+1} 3_{l-1}]) & = & [2_l3_{l+1}] \otimes [3_{l+1}3_{l-1}2_l] + [3_{l-1}2_l] \otimes [2_l3_{l+1}3_{l-1}] \\
&& - \sqrt{[3]}[3_l3_{l+1}] \otimes [3_{l+1}3_{l-1}3_l] - \sqrt{[3]}[3_{l-1}2_l] \otimes [2_l3_{l+1}3_{l-1}].
\end{eqnarray*}
Then from $\mathrm{Im}(\mu_2')_3$ we obtain the following relations in $Y_3$: $[2_0 2_1] \otimes [2_1 2_2 2_0] = [2_{l-1} 2_l] \otimes [2_l 2_{l+1} 2_{l-1}] = -[3_{l-1} 3_l] \otimes [3_l 3_{l+1} 3_{l-1}] = [3_{l-1} 2_l] \otimes [2_l 2_{l+1} 3_{l-1}] = -[2_{l-1} 3_l] \otimes [3_l 3_{l+1} 2_{l-1}] = -[3_{l-1} 2_l] \otimes [2_l 3_{l+1} 3_{l-1}] = [2_{l-1} 3_l] \otimes [3_l 2_{l+1} 2_{l-1}]$, $l=1,2,3$, and thus $Y_3 = \mathbb{C} [2_0 2_1] \otimes [2_1 2_2 2_0] \cong C$.

We now determine $Y_6$. A basis for $(\widetilde{V} \otimes_S A)^S_6$ is given by $[2_{l+1} 2_l] \otimes [2_l 3_{l+1} 2_{l-1} 2_l 2_{l+1}]$, $[3_{l+1} 3_l] \otimes [3_l 2_{l+1} 3_{l-1} 3_l 3_{l+1}]$, $[2_{l+1} 3_l] \otimes [3_l 2_{l+1} 3_{l-1} 3_l 2_{l+1}]$ and $[3_{l+1} 2_l] \otimes [2_l 3_{l+1} 2_{l-1} 2_l 3_{l+1}]$, $l=1,2,3$. For $l=1,2,3$, we have (up to some scalar)
\begin{eqnarray*}
\mu_2'([2_{l+1} 2_l] \otimes [2_l 3_{l+1} 2_{l-1} 2_l 2_{l+1}]) & = & [2_{l-1} 2_l] \otimes [2_l 3_{l+1} 2_{l-1} 2_l 3_{l+1} 2_{l-1}] \\
&& + [2_{l+1} 2_{l-1}] \otimes [2_{l-1} 3_l 2_{l+1} 2_{l-1} 3_l 2_{l+1}], \\
\mu_2'([3_{l+1} 3_l] \otimes [3_l 2_{l+1} 3_{l-1} 3_l 3_{l+1}]) & = & [3_{l-1} 3_l] \otimes [3_l 2_{l+1} 3_{l-1} 3_l 2_{l+1} 3_{l-1}] \\
&& + [3_{l+1} 3_{l-1}] \otimes [3_{l-1} 2_l 3_{l+1} 3_{l-1} 2_l 3_{l+1}], \\
\mu_2'([2_{l+1} 3_l] \otimes [3_l 2_{l+1} 3_{l-1} 3_l 2_{l+1}]) & = & [3_{l-1} 3_l] \otimes [3_l 2_{l+1} 3_{l-1} 3_l 2_{l+1} 3_{l-1}] \\
&& - [2_{l+1} 2_{l-1}] \otimes [2_{l-1} 3_l 2_{l+1} 2_{l-1} 3_l 2_{l+1}], \\
\mu_2'([3_{l+1} 2_l] \otimes [2_l 3_{l+1} 2_{l-1} 2_l 3_{l+1}]) & = & [2_{l-1} 2_l] \otimes [2_l 3_{l+1} 2_{l-1} 2_l 3_{l+1} 2_{l-1}] \\
&& - [3_{l+1} 3_{l-1}] \otimes [3_{l-1} 2_l 3_{l+1} 3_{l-1} 2_l 3_{l+1}],
\end{eqnarray*}
which yield $\mathrm{Im}(\mu_2')_6 = (\widetilde{V} \otimes_S A)^S_6$. Thus $Y_6 = 0$, and we obtain $\overline{HH}_1(A) \cong C$. Then $X_1 = 0$ and from $\chi_{\overline{HC}(A)}(t)$ we deduce that $X_2$ is a graded vector space with Hilbert series $H_{X_2}(t) = t^3 + t^6$.

We now consider $Y' = \overline{HH}_4(A)$, which lives in degrees 8 to 13. Now $\mathrm{ker}(\mu_4') = \mathcal{N}^S[8]$, since $\sum_j w_j^{\ast} \beta(x) w_j = 0$ for all $x \in \mathcal{N}^S_+$ and $\mathcal{N}^S_0 = 0$. As with $\overline{HH}_1(A)$, $Y'_j = 0$ for $j=8,10,11,13$, due to the three-colourability of $\mathcal{E}^{(8)}$. We now determine $Y'_9$. A basis for $(V \otimes_S \mathcal{N})^S_1$ is given by $[2_l 2_{l+1}] \otimes [2_{l+1}]$ and $[3_l 3_{l+1}] \otimes [3_{l+1}]$, $l=1,2,3$, and a basis for $\mathcal{N}^S_1$ is given by $[2_l 2_{l+1}]$ and $[3_l 3_{l+1}]$, $l=1,2,3$. Now $\mu_5'([2_l 2_{l+1}] \otimes [2_{l+1}]) = [2_{l+1} 2_{l-1}] - [2_l 2_{l+1}]$ and $\mu_5'([3_l 3_{l+1}] \otimes [3_{l+1}]) = [3_{l+1} 3_{l-1}] - [3_l 3_{l+1}]$, $l=1,2,3$, thus $Y'_9 = (\mathbb{C}[2_0 2_1] \oplus \mathbb{C}[3_0 3_1])[8]$. We now determine $Y'_{12}$. A basis for $\mathcal{N}^S_4$ is given by $[2_l 3_{l+1} 2_{l-1} 2_l 2_{l+1}]$ and $[3_l 2_{l+1} 3_{l-1} 3_l 3_{l+1}]$, $l=1,2,3$. Since $\mu_5'([1_l 2_{l+1}] \otimes [2_{l+1} 3_{l-1} 3_l 1_{l+1}] = [2_{l+1} 3_{l-1} 2_l 2_{l+1} 2_{l-1}]$ up to some scalar, by using the relations in $A$, and similarly $\mu_5'([4_l 3_{l+1}] \otimes [3_{l+1} 2_{l-1} 2_l 4_{l+1}] = [3_{l+1} 2_{l-1} 3_l 3_{l+1} 3_{l-1}]$, $l=1,2,3$ we see that $Y'_{12} = 0$. Thus $\overline{HH}_4(A) = (\mathbb{C}[2_0 2_1] \oplus \mathbb{C}[3_0 3_1])[8]$, and we obtain $X_3 = (\mathbb{C}[2_0 2_1] \oplus \mathbb{C}[3_0 3_1])[8]$ and we deduce from $\chi_{\overline{HC}(A)}(t)$ that $X_4 = 0$.

To summarize:

\begin{Thm}
The Hochschild and cyclic homology of $A=A(\mathcal{E}^{(8)},W)$, where $W$ is any cell system on $\mathcal{E}^{(8)}$, is given by
$$\begin{array}{lcl}
HH_0(A) \cong S \oplus C, & \qquad & HC_0(A) \cong S \oplus C, \\
HH_1(A) \cong C, && HC_1(A) = 0, \\
HH_2(A) \cong X_2, && HC_2(A) \cong X_2, \\
HH_3(A) \cong X_2, && HC_3(A) = 0, \\
HH_4(A) \cong X_3, && HC_4(A) \cong X_3, \\
HH_5(A) \cong X_3, && HC_5(A) = 0, \\
HH_6(A) \cong X_3^{\ast}[24], && HC_6(A) \cong X_3^{\ast}[24], \\
HH_7(A) \cong X_3^{\ast}[24], && HC_7(A) = 0, \\
HH_8(A) \cong X_2^{\ast}[24], && HC_8(A) \cong X_2^{\ast}[24], \\
HH_9(A) \cong X_2^{\ast}[24], && HC_9(A) = 0, \\
HH_{10}(A) \cong C^{\ast}[24], && HC_{10}(A) \cong C^{\ast}[24], \\
HH_{11}(A) \cong C^{\ast}[24] \oplus K[24], && HC_{11}(A) \cong K[24], \\
HH_{12}(A) \cong C[24] \oplus K^{\ast}[24], && HC_{12}(A) \cong C[24], \\
HH_{12+i}(A) \cong HH_i(A)[24], \quad i \geq 1, && HC_{12+i}(A) \cong HC_i(A)[24], \quad i \geq 1,
\end{array}$$
where the graded vector spaces $C$, $X_2$, $X_3$ and $K$ have Hilbert series $H_C(t) = t^3$, $H_{X_2}(t) = t^3 + t^6$, $H_{X_3}(t) = 2t^9$ and $H_K(t)=2$ respectively.
\end{Thm}

\section{The Hochschild cohomology of $A(\mathcal{G},W)$} \label{sect:Hoch_cohom}
\subsection{The Hochschild cohomology complex} \label{sect:Hoch_cohom-complex}

In this section we will construct a complex which determines the Hochschild cohomology of the almost Calabi-Yau algebra $A = A(\mathcal{G},W)$.
Each four-term piece of this complex will be identified up to a shift in degree with a four-term piece in the Hochschild homology complex (\ref{HH_hom_complex-almostCY}).

The Hochschild cohomology $HH^{\bullet}(A)$ of $A$ may be defined as the derived functor $HH^n(A) = \mathrm{Ext}^n_{A^e}(A,A)$, that is, the homology of the complex
$$0 \rightarrow \mathrm{Hom}_{A^e}(P_0,A) \rightarrow \mathrm{Hom}_{A^e}(P_1,A) \rightarrow \mathrm{Hom}_{A^e}(P_2,A) \rightarrow \cdots$$
where $\quad \cdots \rightarrow P_2 \rightarrow P_1 \rightarrow P_0 \rightarrow A \rightarrow 0$ is any projective resolution of $A$.

Following \cite{etingof/eu:2007}, we can make identifications $\mathrm{Hom}_{A^e}(A \otimes_S \mathcal{N}^{(k)},A) = (\mathcal{N}^{(-k)})^S$, $k=0,1,2$, by identifying $\phi \in \mathrm{Hom}_{A^e}(A \otimes_S \mathcal{N}^{(k)},A)$ with the image $\phi(1 \otimes 1) = x \in (\mathcal{N}^{(-k)})^S$. We write $\phi = x \circ - : A \otimes_S \mathcal{N}^{(-k)} \rightarrow A$, and have $\phi(y \otimes z) = x \circ (y \otimes z) = yx\beta^{k}(z)$, for $x \in (\mathcal{N}^{(-k)})^S$, $y \in A$, $z \in \mathcal{N}^{(k)}$.
We also make identifications $\mathrm{Hom}_{A^e}(A \otimes_S V \otimes_S \mathcal{N}^{(k)},A) = (\widetilde{V} \otimes_S \mathcal{N}^{(-k)})^S[-2]$, $k=0,1,2$, by identifying $\phi \in \mathrm{Hom}_{A^e}(A \otimes_S V \otimes_S \mathcal{N}^{(k)},A)$ which maps $1 \otimes a \otimes 1 \mapsto x_a$ with the element $\sum_{a \in \mathcal{G}_1} \widetilde{a} \otimes x_a \in (\widetilde{V} \otimes_S \mathcal{N}^{(-k)})^S$. We write $\phi = \sum_{b \in \mathcal{G}_1} \widetilde{b} \otimes x_b \circ - : A \otimes_S V \otimes_S \mathcal{N}^{(-k)} \rightarrow A$, and have $\phi(y \otimes a \otimes z) = \sum_{b \in \mathcal{G}_1} \widetilde{b} \otimes x_b \circ (y \otimes a \otimes z) = yx_a\beta^{k}(z)$, for $\widetilde{a} \otimes x_a \in (\widetilde{V} \otimes_S \mathcal{N}^{(-k)})^S$, $y \in A$, $z \in \mathcal{N}^{(k)}$.
Similarly, we identify $\mathrm{Hom}_{A^e}(A \otimes_S \widetilde{V} \otimes_S \mathcal{N}^{(k)},A) = (V \otimes_S \mathcal{N}^{(-k)})^S[-2]$, $k=0,1,2$, by identifying $\phi$ which maps $1 \otimes \widetilde{a} \otimes 1 \mapsto y_a$ with the element $\sum_{a \in \mathcal{G}_1} a \otimes y_a$. We write $\phi = \sum_{b \in \mathcal{G}_1} b \otimes y_b \circ - : A \otimes_S \widetilde{V} \otimes_S \mathcal{N}^{(-k)} \rightarrow A$, and have $\phi(y \otimes \widetilde{a} \otimes z) = yy_a\beta^{k}(z)$, for $a \otimes y_a \in (V \otimes_S \mathcal{N}^{(-k)})^S$, $y \in A$, $z \in \mathcal{N}^{(k)}$.

Applying the functor $\mathrm{Hom}_{A^e}(-,A)$ to the periodic resolution (\ref{resolution-almostCY}) we get the Hochschild cohomology complex:
\begin{eqnarray}
&& (\mathcal{N}^{(2)})^S[-h] \stackrel{\mu_{4}^{\ast}}{\leftarrow} A^S[-3] \stackrel{\mu_{3}^{\ast}}{\leftarrow} (V \otimes_S A)^S[-3] \stackrel{\mu_{2}^{\ast}}{\leftarrow} (\widetilde{V} \otimes_S A)^S[-2] \stackrel{\mu_{1}^{\ast}}{\leftarrow} A^S \leftarrow 0 \nonumber \\
&& \mathcal{N}^S[-2h] \stackrel{\mu_{8}^{\ast}}{\leftarrow} (\mathcal{N}^{(2)})^S[-h-3] \stackrel{\mu_{7}^{\ast}}{\leftarrow} (V \otimes_S \mathcal{N}^{(2)})^S[-h-3] \stackrel{\mu_{6}^{\ast}}{\leftarrow} (\widetilde{V} \otimes_S \mathcal{N}^{(2)})^S[-h-2] \stackrel{\mu_{5}^{\ast}}{\leftarrow} \nonumber \\
& \cdots & \leftarrow A^S[-3h] \stackrel{\mu_{12}^{\ast}}{\leftarrow} \mathcal{N}^S[-2h-3] \stackrel{\mu_{11}^{\ast}}{\leftarrow} (V \otimes_S \mathcal{N})^S[-2h-3] \stackrel{\mu_{10}^{\ast}}{\leftarrow} (\widetilde{V} \otimes_S \mathcal{N})^S[-2h-2] \stackrel{\mu_{9}^{\ast}}{\leftarrow} \label{HH_cohom_complex-almostCY}
\end{eqnarray}

\begin{Prop}
We have $\mu_i^{\ast} = \pm\mu_{16-i}'$.
\end{Prop}

\noindent \emph{Proof}:
(i) $\mu_1^{\ast} = -\mu_3'$:
Let $a \in V$ and $x \in A^S$. Then
$$\mu_1^{\ast}(x)(1 \otimes a \otimes 1) = x \circ \mu_1(1 \otimes a \otimes 1) = x \circ (a \otimes 1 - 1 \otimes a) = ax-xa.$$
So $\mu_1^{\ast}(x)$ maps $1 \otimes a \otimes 1 \mapsto [a,x]$, giving $\mu_1^{\ast}(x) = \sum_{a \in \mathcal{G}_1} \widetilde{a} \otimes [a,x] = -\mu_3'(x)$.
Similarly, $\mu_5^{\ast}(x)$ maps $1 \otimes a \otimes 1 \mapsto ax-x\beta(a)$, giving $\mu_5^{\ast}(x) = \sum_{a \in \mathcal{G}_1} \widetilde{a} \otimes (ax-x\beta(a)) = -\mu_{11}'(x)$, and we also have $\mu_9^{\ast}(x) = \sum_{a \in \mathcal{G}_1} \widetilde{a} \otimes (ax-x\beta^2(a)) = -\mu_{7}'(x)$. \\
(ii) $\mu_2^{\ast} = \mu_2'$:
Let $a' \in V$ and for each $a \in V$ let $x_a$ be a homogeneous element in $A$ such that $\widetilde{a} \otimes x_a \in (\widetilde{V} \otimes A)^S$. Then
\begin{eqnarray*}
\lefteqn{\mu_2^{\ast}(\sum_{a \in \mathcal{G}_1} \widetilde{a} \otimes x_a)(1 \otimes \widetilde{a'} \otimes 1) \;\; = \;\; \sum_{a \in \mathcal{G}_1} \widetilde{a} \otimes x_a \circ \mu_2(1 \otimes \widetilde{a'} \otimes 1)} \\
& = & \sum_{a \in \mathcal{G}_1} \widetilde{a} \otimes x_a \circ \left( \sum_{b,b' \in \mathcal{G}_1} W_{a'bb'} (b \otimes b' \otimes 1 + 1 \otimes b \otimes b') \right) \;\; = \;\; \sum_{b,b' \in \mathcal{G}_1} W_{a'bb'} (bx_{b'}+x_{b}b').
\end{eqnarray*}
So $\mu_2^{\ast}(\sum_{a \in \mathcal{G}_1} \widetilde{a} \otimes x_a)$ maps $1 \otimes \widetilde{a'} \otimes 1 \mapsto \sum_{b,b' \in \mathcal{G}_1} W_{a'bb'} (bx_{b'}+x_{b}b')$, giving $\mu_2^{\ast}(\sum_{a \in \mathcal{G}_1} \widetilde{a} \otimes x_a) = \sum_{a,b,b' \in \mathcal{G}_1} W_{abb'} (a \otimes bx_{b'} + a \otimes x_{b}b') = \mu_2'(\sum_{a \in \mathcal{G}_1} \widetilde{a} \otimes x_a)$.
Similarly, $\mu_6^{\ast}(\sum_{a \in \mathcal{G}_1} \widetilde{a} \otimes x_a) = \sum_{a,b,b' \in \mathcal{G}_1} W_{abb'} (a \otimes bx_{b'} + a \otimes x_{b}\beta(b')) = \mu_{10}'(\sum_{a \in \mathcal{G}_1} \widetilde{a} \otimes x_a)$ and $\mu_{10}^{\ast}(\sum_{a \in \mathcal{G}_1} \widetilde{a} \otimes x_a) = \sum_{a,b,b' \in \mathcal{G}_1} W_{abb'} (a \otimes bx_{b'} + a \otimes x_{b}\beta^2(b')) = \mu_{6}'(\sum_{a \in \mathcal{G}_1} \widetilde{a} \otimes x_a)$. \\
(iii) $\mu_3^{\ast} = -\mu_1'$:
For each $a \in V$ let $y_a$ be a homogeneous element in $A$ such that $a \otimes y_a \in (V \otimes A)^S$. Then
\begin{eqnarray*}
\mu_3^{\ast}(\sum_{a \in \mathcal{G}_1} a \otimes y_a)(1 \otimes 1) & = & \sum_{a \in \mathcal{G}_1} a \otimes y_a \circ \mu_3(1 \otimes 1) \\
& = & \sum_{a \in \mathcal{G}_1} a \otimes y_a \circ \sum_{b \in \mathcal{G}_1} (b \otimes \widetilde{b} \otimes 1 - 1 \otimes \widetilde{b} \otimes b) \;\; = \;\; \sum_{b \in \mathcal{G}_1} (by_{b}-y_{b}b).
\end{eqnarray*}
So $\mu_3^{\ast}(\sum_{a \in \mathcal{G}_1} a \otimes y_a)$ maps $1 \otimes 1 \mapsto \sum_{b \in \mathcal{G}_1} [b,y_b]$, giving $\mu_3^{\ast}(\sum_{a \in \mathcal{G}_1} a \otimes y_a) = \sum_{a \in \mathcal{G}_1} [a,y_a] = -\mu_1'(\sum_{a \in \mathcal{G}_1} a \otimes y_a)$.
Similarly, $\mu_7^{\ast}(\sum_{a \in \mathcal{G}_1} a \otimes y_a) = \sum_{a \in \mathcal{G}_1} (ay_{a}-y_{a}\beta(a)) = -\mu_9'(\sum_{a \in \mathcal{G}_1} a \otimes y_a)$ and $\mu_{11}^{\ast}(\sum_{a \in \mathcal{G}_1} a \otimes y_a) = \sum_{a \in \mathcal{G}_1} (ay_{a}-y_{a}\beta^2(a)) = -\mu_5'(\sum_{a \in \mathcal{G}_1} a \otimes y_a)$. \\
(iv) $\mu_4^{\ast} = \mu_{12}'$:
Let $x \in A^S$. Then
$$\mu_4^{\ast}(x)(1 \otimes 1) = x \circ \mu_4(1 \otimes 1) = x \circ \sum_j w_j \otimes w_j^{\ast} = \sum_j w_j x w_j^{\ast},$$
where $\{ w_j \}$ is a homogeneous basis for $A$ and $\{ w_j^{\ast} \}$ is its corresponding dual basis.
So $\mu_4^{\ast}(x)$ maps $1 \otimes 1 \mapsto \sum_j w_j x w_j^{\ast}$, giving $\mu_4^{\ast}(x) = \sum_j w_j x w_j^{\ast} = \mu_{12}'(x)$. Similarly, $\mu_8^{\ast}(x) = \sum_j w_j x \beta(w_j^{\ast}) = \mu_{8}'(x)$ and $\mu_{12}^{\ast}(x) = \sum_j w_j x \beta^2(w_j^{\ast}) = \mu_4'(x)$.
\hfill
$\Box$

Thus we see that we can identify, up to a shift in degree, each four-term portion of the cohomology complex (\ref{HH_cohom_complex-almostCY}) with a portion of the homology complex (\ref{HH_hom_complex-almostCY}):
\begin{eqnarray*}
HH^i(A) & \cong & HH_{3-i}(A)[-3], \qquad \qquad \quad i=1,2, \\
HH^i(A) & \cong & HH_{15-i}(A)[-3h-3], \qquad \; i=3,\ldots,12, \\
HH^{12+i}(A) & \cong & HH^i(A)[-3h], \qquad \qquad \quad \; i=1,2,\ldots,
\end{eqnarray*}
and the self-duality of the homology complex (\ref{HH_hom_complex-almostCY}) yields the relations
\begin{eqnarray*}
HH^i(A)^{\ast} & \cong & HH^{7-i}(A), \qquad \qquad i=1,\ldots,6, \\
HH^i(A)^{\ast} & \cong & HH^{19-i}(A), \qquad \qquad i=7,\ldots,11.
\end{eqnarray*}

\subsection{The Hochschild cohomology of $A = A(\mathcal{G},W)$} \label{sect:Hoch_cohom-results}

For $HH^0(A) = \mathrm{ker}(\mu_1^{\ast})/\mathrm{Im}(\mu_0^{\ast}) = \mathrm{ker}(\mu_1^{\ast})$, we have $HH^0(A) \cong HH_3(A)'[-3] \oplus L$, where $HH_3(A)' = \oplus_{j=3}^{h-1} HH_3(A)_j$ and $L = \mathbb{C} \{ u_{j\nu(j)} | \, \nu(j) = j \}$.

Then we have the following results for the Hochschild cohomology of $A$:

\begin{Thm}
The Hochschild cohomology of $A=A(\mathcal{A}^{(4)},W)$, where $W$ is any cell system on $\mathcal{A}^{(4)}$, is given by
$$\begin{array}{lcl}
HH^0(A) \cong X[-3], & \qquad & HH^1(A) \cong X[-3], \\
HH^6(A) \cong X^{\ast}[-3], && HH^7(A) \cong X^{\ast}[-3], \\
HH^{12}(A) \cong X[-15], && HH^{j}(A) = 0, \quad j = 2,\ldots,5,8,\ldots,11,
\end{array}$$
and $HH^{12+i}(A) \cong HH^i(A)[-12]$ for $i \geq 1$, where the graded vector space $X$ has Hilbert series $H_X(t) = t^3$.
\end{Thm}

\begin{Thm}
The Hochschild cohomology of $A=A(\mathcal{A}^{(5)},W)$, where $W$ is any cell system on $\mathcal{A}^{(5)}$, is given by
$$\begin{array}{lcl}
HH^0(A) \cong X_2[-3], & \qquad & HH^1(A) \cong X_2[-3], \\
HH^2(A) = 0, && HH^3(A) = 0, \\
HH^4(A) = 0, && HH^5(A) = 0, \\
HH^6(A) \cong X_2^{\ast}[-3], && HH^7(A) \cong X_2^{\ast}[-3], \\
HH^8(A) \cong X_3^{\ast}[-3], && HH^9(A) \cong X_3^{\ast}[-3], \\
HH^{10}(A) \cong X_3[-18], && HH^{11}(A) \cong X_3[-18], \\
HH^{12}(A) \cong X_2[-18], && HH^{12+i}(A) \cong HH^i(A)[-15], \quad i \geq 1,
\end{array}$$
where the graded vector spaces $X_2$ and $X_3$ have Hilbert series $H_{X_2}(t) = t^3$ and $H_{X_3}(t) = t^6$ respectively.
\end{Thm}

\begin{Thm}
The Hochschild cohomology of $A=A(\mathcal{A}^{(6)},W)$, where $W$ is any cell system on $\mathcal{A}^{(6)}$, is given by
$$\begin{array}{lcl}
HH^0(A) \cong X[-3] \oplus L, & \qquad & HH^1(A) \cong X[-3], \\
HH^2(A) = 0, && HH^3(A) \cong K^{\ast}[-3], \\
HH^4(A) \cong K[-3], && HH^5(A) = 0, \\
HH^6(A) \cong X^{\ast}[-3], && HH^7(A) \cong X^{\ast}[-3], \\
HH^8(A) = 0, && HH^9(A) = 0, \\
HH^{10}(A) = 0, && HH^{11}(A) = 0, \\
HH^{12}(A) \cong X[-21], && HH^{12+i}(A) \cong HH^i(A)[-18], \quad i \geq 1,
\end{array}$$
where the graded vector spaces $L$, $X$ and $K$ have Hilbert series $H_L(t) = t^3$, $H_{X}(t) = t^3$ and $H_K(t)=2$ respectively.
\end{Thm}

\begin{Thm}
The Hochschild cohomology of $A=A(\mathcal{A}^{(7)},W)$, where $W$ is any cell system on $\mathcal{A}^{(7)}$, is given by
$$\begin{array}{lcl}
HH^0(A) \cong X_2[-3], & \qquad & HH^1(A) \cong X_2[-3], \\
HH^2(A) = 0, && HH^3(A) \cong K^{\ast}[-3], \\
HH^4(A) \cong K[-3], && HH^5(A) = 0, \\
HH^6(A) \cong X_2^{\ast}[-3], && HH^7(A) \cong X_2^{\ast}[-3], \\
HH^8(A) \cong X_3^{\ast}[-3], && HH^9(A) \cong X_3^{\ast}[-3], \\
HH^{10}(A) \cong X_3[-24], && HH^{11}(A) \cong X_3[-24], \\
HH^{12}(A) \cong X_2[-24], && HH^{12+i}(A) \cong HH^i(A)[-21], \quad i \geq 1,
\end{array}$$
where the graded vector spaces $X_2$, $X_3$ and $K$ have Hilbert series $H_{X_2}(t) = t^3 + t^6$, $H_{X_3}(t) = t^9$ and $H_K(t)=2$ respectively.
\end{Thm}

\begin{Thm}
The Hochschild cohomology of $A=A(\mathcal{D}^{(6k)},W)$, $k \geq 1$, where $W$ is equivalent to one of the cell systems given in \cite{evans/pugh:2009i},  is given by
$$\begin{array}{lcl}
HH^0(A) \cong C^{\ast}[6k-3] \oplus L, & \qquad & HH^1(A) \cong C^{\ast}[6k-3] \oplus X^{\ast}[6k-3], \\
HH^2(A) \cong C[-3] \oplus X[-3], && HH^3(A) \cong C[-3] \oplus K^{\ast}[-3], \\
HH^4(A) \cong C^{\ast}[-3] \oplus K[-3], && HH^{4+i}(A) \cong HH^i(A)[-6k], \quad i \geq 1,
\end{array}$$
where the graded vector spaces $C$, $L$, $X$, $K$ have Hilbert series $H_C(t) = \sum_{j=1}^{2k-2} 3t^{3j} + t^{6k-3}$, $H_L(t) = (3k(2k-1)+3)t^{6k-3}$, $H_X(t) = t^3 + \sum_{j=2}^{2k-2} 3t^{3j} + t^{3k} + t^{6k-3}$ and $H_K(t)=6k(k-1)+2$ respectively, where for $k=1$, $H_X(t) = 0$.
\end{Thm}

\begin{Thm}
The Hochschild cohomology of $A=A(\mathcal{D}^{(6k+3)},W)$, $k \geq 1$, where $W$ is equivalent to one of the cell systems given in \cite{evans/pugh:2009i},  is given by
$$\begin{array}{lcl}
HH^0(A) \cong C^{\ast}[6k] \oplus L, & \qquad & HH^1(A) \cong C^{\ast}[6k] \oplus X^{\ast}[6k], \\
HH^2(A) \cong C[-3] \oplus X[-3], && HH^3(A) \cong C[-3] \oplus K^{\ast}[-3], \\
HH^4(A) \cong C^{\ast}[-3] \oplus K[-3], && HH^{4+i}(A) \cong HH^i(A)[-6k-3], \quad i \geq 1,
\end{array}$$
where the graded vector spaces $C$, $L$, $X$, $K$ have Hilbert series $H_C(t) = \sum_{j=1}^{2k-1} 3t^{3j} + t^{6k}$, $H_L(t) = (3k(2k+1)+3)t^{6k}$, $H_X(t) = t^3 + \sum_{j=2}^{2k-1} 3t^{3j} + t^{6k}$ and $H_K(t)=6k^2$ respectively.
\end{Thm}

\begin{Thm}
The Hochschild cohomology of $A=A(\mathcal{A}^{(n)\ast},W)$, $n \geq 4$, where $W$ is any cell system on $\mathcal{A}^{(n)\ast}$, is given by
$$\begin{array}{lcl}
HH^0(A) \cong C^{\ast}[n-3] \oplus L, & \qquad & HH^1(A) \cong C^{\ast}[n-3], \\
HH^2(A) \cong C[-3], && HH^3(A) \cong C[-3], \\
HH^4(A) \cong C^{\ast}[-3], && HH^{4+i}(A) \cong HH^i(A)[-n], \quad i \geq 1,
\end{array}$$
where the graded vector spaces $C$, $L$ have Hilbert series $H_C(t) = \sum_{j=1}^{n-3} \lfloor (n-j-1)/2 \rfloor t^{j}$ and $H_L(t) = \lfloor (n-1)/2 \rfloor t^{n-3}$.
\end{Thm}

\begin{Thm}
The Hochschild cohomology of $A=A(\mathcal{D}^{(6k)\ast},W)$, $k \geq 1$, where $W$ is equivalent to one of the cell systems given in \cite{evans/pugh:2009i},  is given by
$$\begin{array}{lcl}
HH^0(A) \cong C^{\ast}[6k-3] \oplus L, & \qquad & HH^1(A) \cong C^{\ast}[6k-3], \\
HH^2(A) \cong C[-3], && HH^3(A) \cong C[-3] \oplus K^{\ast}[-3], \\
HH^4(A) \cong C^{\ast}[-3] \oplus K[-3], && HH^{4+i}(A) \cong HH^i(A)[-6k], \quad i \geq 1,
\end{array}$$
where the graded vector spaces $C$, $L$, $K$ have Hilbert series $H_C(t) = \sum_{j=1}^{\lfloor (2m-1)/3 \rfloor} (m- \lfloor 3j/2 \rfloor) t^{3j}$, $H_L(t) = (9k-3)t^{6k-3}$ and $H_K(t)=6k-4$ respectively.
\end{Thm}

\begin{Thm}
The Hochschild cohomology of $A=A(\mathcal{D}^{(6k+3)},W)$, $k \geq 1$, where $W$ is equivalent to one of the cell systems given in \cite{evans/pugh:2009i}, is given by
$$\begin{array}{lcl}
HH^0(A) \cong C^{\ast}[6k] \oplus L, & \qquad & HH^1(A) \cong C^{\ast}[6k], \\
HH^2(A) \cong C[-3], && HH^3(A) \cong C[-3] \oplus K^{\ast}[-3], \\
HH^4(A) \cong C^{\ast}[-3] \oplus K[-3], && HH^{4+i}(A) \cong HH^i(A)[-6k-3], \quad i \geq 1,
\end{array}$$
where the graded vector spaces $C$, $L$, $K$ have Hilbert series $H_C(t) = \sum_{j=1}^{\lfloor (2m-2)/3 \rfloor} (m- \lfloor (3j+1)/2 \rfloor) t^{3j}$, $H_L(t) = (9k+3)t^{6k}$ and $H_K(t)=6k$ respectively.
\end{Thm}

\begin{Thm}
The Hochschild cohomology of $A=A(\mathcal{E}^{(8)},W)$, where $W$ is any cell system on $\mathcal{E}^{(8)}$, is given by
$$\begin{array}{lcl}
HH^0(A) \cong X_2[-3], & \qquad & HH^1(A) \cong X_2[-3], \\
HH^2(A) \cong C[-3], && HH^3(A) \cong C[-3] \oplus K^{\ast}[-3], \\
HH^4(A) \cong C^{\ast}[-3] \oplus K[-3], && HH^5(A) \cong C^{\ast}[-3], \\
HH^6(A) \cong X_2^{\ast}[-3], && HH^7(A) \cong X_2^{\ast}[-3], \\
HH^8(A) \cong X_3^{\ast}[-3], && HH^9(A) \cong X_3^{\ast}[-3], \\
HH^{10}(A) \cong X_3[-27], && HH^{11}(A) \cong X_3[-27], \\
HH^{12}(A) \cong X_2[-27], && HH^{12+i}(A) \cong HH^i(A)[-24], \quad i \geq 1,
\end{array}$$
where the graded vector spaces $C$, $X_2$, $X_3$ and $K$ have Hilbert series $H_C(t) = t^3$, $H_{X_2}(t) = t^3 + t^6$, $H_{X_3}(t) = 2t^9$ and $H_K(t)=2$ respectively.
\end{Thm}

\begin{Thm}
The Hochschild cohomology of $A=A(\mathcal{E}^{(8)\ast},W)$, where $W$ is any cell system on $\mathcal{E}^{(8)\ast}$, is given by
$$\begin{array}{lcl}
HH^0(A) \cong C^{\ast}[5] \oplus L, & \qquad & HH^1(A) \cong C^{\ast}[5], \\
HH^2(A) \cong C[-3], && HH^3(A) \cong C[-3] \oplus K^{\ast}[-3], \\
HH^4(A) \cong C^{\ast}[-3] \oplus K[-3], && HH^{4+i}(A) \cong HH^i(A)[-8], \quad i \geq 1,
\end{array}$$
where the graded vector spaces $C$, $L$, $K$ have Hilbert series $H_C(t) = 2t + t^2 + t^3 + t^5$, $H_L(t) = 4t^5$ and $H_K(t)=2$ respectively.
\end{Thm}

\paragraph{Acknowledgements}

Both authors were supported by the Marie Curie Research Training Network MRTN-CT-2006-031962 EU-NCG.
The authors would like to thank Karin Erdmann, Pavel Etingof, Victor Ginzburg and Jean-Louis Loday for helpful discussions and correspondence.

\end{document}